\documentclass[10pt,a4paper]{article}
\usepackage{amsmath,amsthm,amssymb,latexsym,epic}
\usepackage{bezier}

\newtheorem{theorem}{Theorem}
\newtheorem{lemma}[theorem]{Lemma}
\newtheorem{corollary}[theorem]{Corollary}
\newtheorem{proposition}[theorem]{Proposition}

\usepackage[all]{xy}
\usepackage{enumerate}
\newcommand{\ls}{\mathrm{ls}}

\newcommand{\GR}{\mathcal{R}}
\newcommand{\GL}{\mathcal{L}}
\newcommand{\GH}{\mathcal{H}}
\newcommand{\GJ}{\mathcal{J}}
\newcommand{\GD}{\mathcal{D}}
\newcommand{\IS}{\mathcal{IS}}

\newcommand{\T}{\mathcal{T}}
\newcommand{\Sym}{\mathcal{S}}

\newcommand{\corank}{\mathrm{corank}}
\newcommand{\lmod}{\mid\!}
\newcommand{\rmod}{\!\mid}
\newcommand{\B}{\mathfrak{B}}
\newcommand{\A}{\mathcal{A}}
\newcommand{\len}{\mathfrak{l}}
\newcommand{\C}{\mathfrak{C}}
\newcommand{\PB}{\mathcal{P}\mathfrak{B}}
\newcommand{\bn}{{\bf n}}

\title{Presentation of the singular part of the Brauer monoid}
\author{Victor Maltcev and Volodymyr Mazorchuk}

\date{}

\begin{document}

\maketitle

\begin{abstract}
We obtain a presentation for the singular part of the Brauer monoid
with respect to an irreducible system of generators, consisting of 
idempotents. As an application of this result we get a new 
construction of the symmetric group via connected sequences of
subsets. Another application describes the lengths of elements in
the singular part of the Brauer monoid with respect to the system of 
generators, mentioned above.
\end{abstract}

\section{Introduction}\label{sec:intro}

The symmetric group $\Sym_n$ is a central object of study in many 
branches of mathematics. There exist several ``natural'' analogues 
(or generalizations) of  $\Sym_n$ in the theory of semigroups. 
The most classical ones are the symmetric semigroup $\T_n$ and the 
inverse symmetric  semigroup $\IS_n$. These arise when one tries to 
generalize Cayley's Theorem to the classes of all semigroups or all
inverse semigroups respectively. A less obvious semigroup generalization
of $\Sym_n$ is the so-called {\em Brauer semigroup} $\B_n$, which
appears in the context of centralizer algebras in representation theory,
see \cite{Br}. $\B_n$ contains $\Sym_n$ as the subgroup of all 
invertible elements and has a  nice geometric realization (see
Section~\ref{sec:definition}). The deformation of the corresponding
semigroup algebra, the so-called {\em Brauer algebra}, has been
intensively studied by specialists in representation theory, knot
theory and theoretical physics. The semigroup properties of 
$\B_n$ were studied in \cite{Maz1,Maz2,Mal,KM1,KMM,KM2}.

Given a finitely generated semigroup, a fundamental question is to
find its presentation with respect to some (irreducible) system of
generators. For example, for $\Sym_n$ and $\B_n$ several such 
presentations are known. However, for semigroups one can even make
the problem more semigroup-oriented, and ask to find a presentation for
the {\em singular part} of the semigroup, which, by definition, is the set
of all non-invertible elements. In the case of a finite semigroup
all non-invertible elements form again a semigroup and hence the
problem to find a presentation for the singular part makes sense. 
For example, in \cite{East} a presentation for the singular part 
of  $\IS_n$ is found  (a presentation for $\IS_n$ itself can be 
found in \cite{Aiz}).

From \cite{Mal} we know that $\B_n\setminus \Sym_n$ has a natural
irreducible system of generators, consisting of idempotents. The
main aim of the present paper is to obtain a presentation of 
$\B_n\setminus \Sym_n$ with respect to this system of generators.
Surprisingly enough, the system of the corresponding defining relations
is not big and all relations have an obvious interpretation via the
geometric realization of $\B_n$. This result is presented in 
Theorem~\ref{th:main}. As usual, a  tricky part in the proof of 
Theorem~\ref{th:main} is to show that the listed system of 
defining relations is complete.
This part of the proof is quite technical and occupies the whole
Section~\ref{sec:main}. In Section~\ref{sec:application} we present
several combinatorial applications of Theorem~\ref{th:main}. 
These include an interesting combinatorial realization of the 
symmetric group via equivalence classes of sequences of
``connected'' two-element subsets, and a computation of the 
maximal length for an element in $\B_n\setminus \Sym_n$
with respect to our system of generators.

\bigskip
{\bf Acknowledgments.} The paper was written during the visit of the 
first author to Uppsala University, which was supported by the 
Swedish Institute. The financial support of the Swedish Institute and 
the hospitality of Uppsala University are greatfully acknowledged. For 
the second author the research was partially supported by the 
Swedish Research Council.

%%%%%%%%%%%%%%%%%%%%%%%%%%%%%%%%%%%%%%%%%%%%%%%%%%%%%%%%%%%%%%%%%%%%%%%%%%%%%%%%%%%%%%%%%%%%%%%%%%%%

\section{Preliminaries about $\B_n$}\label{sec:definition}

Let $n$ be a positive integer. Put $\bn=\{1,\dots,n\}$ and
$\bn'=\{1',\dots,n'\}$. We consider the map $':\bn\to\bn'$ as a 
fixed bijection and denote the inverse bijection by the same symbol, 
that is $(x')'=x$ for all $x\in\bn$. 
The elements of the Brauer semigroup $\B_n$ are all possible 
partitions of $\bn\cup\bn'$ into two-element blocks. 
It is easy to see that $\lmod\B_n\rmod=(2n-1)!!$.

A two-element subset, $\{i,j\}$, of $\bn\cup\bn'$ will be called a
\begin{itemize}
\item
\emph{left bracket} provided that $\{i,j\}\subset \bn$;
\item
\emph{right bracket} provided that $\{i,j\}\subset \bn'$;
\item
\emph{line}, if $\{i,j\}$ is neither a left nor a right bracket.
\end{itemize}
Obviously, every element of $\B_n$ contains the same number of 
left and right brackets. Let $\pi\in\B_n$, and assume that 
$\{i_k,j_k\}$, $k\in K$, is a list of all left brackets of $\pi$;
$\{u_k',v_k'\}$, $k\in K$, is a list of all right brackets of $\pi$;
and $\{f_l,g_l'\}$, $l\in L$, is a list of all lines of $\pi$.
Then we have 
\begin{equation}\label{eq:convenient-form}
\pi=\bigl\{\{i_k,j_k\}_{k\in K},\{u_k',v_k'\}_{k\in K},\{f_l,g_l'\}_{l\in L}\bigr\}.
\end{equation}
We say that $\pi$ has \emph{corank} $\corank(\pi)=2\lmod K\rmod\leq 
2\lfloor\frac{n}{2}\rfloor$. 

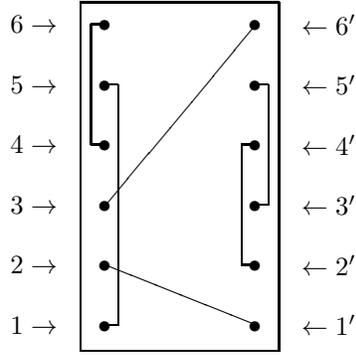
\begin{figure}
%%%%%%%%%%%%%%%%%%%%%%%%%%%%%%%%%%%%%%%%%%%%%%%%%%%%%%%%%%%%%%%%%%%%%%%%%%%%%%%%
%%%%%
%\input{figa1.pic}
%TexCad Options
%\grade{\on}
%\paths{\on}
%\beziermacro{\off}
%\reduce{\on}
%\snapping{\off}
%\quality{2.00}
%\graddiff{0.01}
%\snapasp{1}
%\zoom{1.00}
\special{em:linewidth 0.4pt} \unitlength 0.80mm
\linethickness{1pt}
\begin{picture}(150.00,55.00)
\put(65.00,00.00){\makebox(0,0)[cc]{$\bullet$}}
\put(65.00,10.00){\makebox(0,0)[cc]{$\bullet$}}
\put(65.00,20.00){\makebox(0,0)[cc]{$\bullet$}}
\put(65.00,30.00){\makebox(0,0)[cc]{$\bullet$}}
\put(65.00,40.00){\makebox(0,0)[cc]{$\bullet$}}
\put(65.00,50.00){\makebox(0,0)[cc]{$\bullet$}}
\put(90.00,10.00){\makebox(0,0)[cc]{$\bullet$}}
\put(90.00,20.00){\makebox(0,0)[cc]{$\bullet$}}
\put(90.00,30.00){\makebox(0,0)[cc]{$\bullet$}}
\put(90.00,40.00){\makebox(0,0)[cc]{$\bullet$}}
\put(90.00,50.00){\makebox(0,0)[cc]{$\bullet$}}
\put(90.00,00.00){\makebox(0,0)[cc]{$\bullet$}}
%%%%%%%%%%%%%%%%%%%%%%%%%%%%%%%%%%%%%%%%%%%%%%%%%%%%%%%%% FRAME %%%%%%%%%%%%%%%%%%%%%%%%%%%%%%%%

\drawline(61.00,-04.00)(61.00,54.00)
\drawline(61.00,54.00)(94.00,54.00)
\drawline(94.00,54.00)(94.00,-04.00)
\drawline(94.00,-04.00)(61.00,-04.00)

%%%%%%%%%%%%%%%%%%%%%%%%%%%%%%%%%%%%%%%%%%%%%%%%%%%%%%%% THE ELEMENT %%%%%%%%%%%%%%%%%%%%%%%%%%%%

\drawline(65.30,20.30)(90.30,50.30)
\drawline(65.10,10.30)(90.30,00.30)

\drawline(65.30,00.30)(67.30,00.30)
\drawline(67.30,00.30)(67.30,40.30)
\drawline(67.30,40.30)(65.30,40.30)

\drawline(65.30,30.30)(62.80,30.30)
\drawline(62.80,30.30)(62.80,50.30)
\drawline(62.80,50.30)(65.30,50.30)

\drawline(90.30,10.30)(87.80,10.30)
\drawline(87.80,10.30)(87.80,30.30)
\drawline(87.80,30.30)(90.30,30.30)

\drawline(90.30,20.30)(92.30,20.30)
\drawline(92.30,20.30)(92.30,40.30)
\drawline(92.30,40.30)(90.30,40.30)

%%%%%%%%%%%%%%%%%%%%%%%%%%%%%%%%%%%%%%%%%%%%%%%%%%%%%%%%%%%%%%%%%%%%%%%%%%%%%%%%%%%%%%%%%
\put(53.30,00.30){\makebox(0,0)[cc]{$1\rightarrow$}}
\put(53.30,10.30){\makebox(0,0)[cc]{$2\rightarrow$}}
\put(53.30,20.30){\makebox(0,0)[cc]{$3\rightarrow$}}
\put(53.30,30.30){\makebox(0,0)[cc]{$4\rightarrow$}}
\put(53.30,40.30){\makebox(0,0)[cc]{$5\rightarrow$}}
\put(53.30,50.30){\makebox(0,0)[cc]{$6\rightarrow$}}
%%%%%%%%%%%%%%%%%%%%%%%%%%%%%%%%%%%%%%%%%%%%%%%%%%%%%%%%%%%%%%%%%%%%%%%%%%%%%%%%%%%%%%%
\put(102.30,00.30){\makebox(0,0)[cc]{$\leftarrow 1'$}}
\put(102.30,10.30){\makebox(0,0)[cc]{$\leftarrow 2'$}}
\put(102.30,20.30){\makebox(0,0)[cc]{$\leftarrow 3'$}}
\put(102.30,30.30){\makebox(0,0)[cc]{$\leftarrow 4'$}}
\put(102.30,40.30){\makebox(0,0)[cc]{$\leftarrow 5'$}}
\put(102.30,50.30){\makebox(0,0)[cc]{$\leftarrow 6'$}}
\end{picture}
\caption{The element $\bigl\{\{1,5\},\{4,6\},\{2',4'\},\{3',5'\},\{2,1'\},\{3,6'\}\bigr\}$ of $\B_6$.}\label{fig:example-element}
%%%%%%%%%%%%%%%%%%%%%%%%%%%%%%%%%%%%%%%%%%%%%%%%%%%%%%%%%%%%%%%%%%%%%%%%%%%%%%%%
%%%%
\end{figure}

It is convenient to represent the elements of $\B_n$ geometrically
as a kind of {\em microchips} as follows: we have two sets of pins 
(which correspond to elements
in $\bn$ and $\bn'$ respectively), which are connected in pairs
(this corresponds to the partition of $\bn\cup\bn'$ into 
two-element blocks, which our element from $\B_n$ represents).
An example is shown on Figure~\ref{fig:example-element}, for the
convenience the same element is also written in
the form \eqref{eq:convenient-form}. 

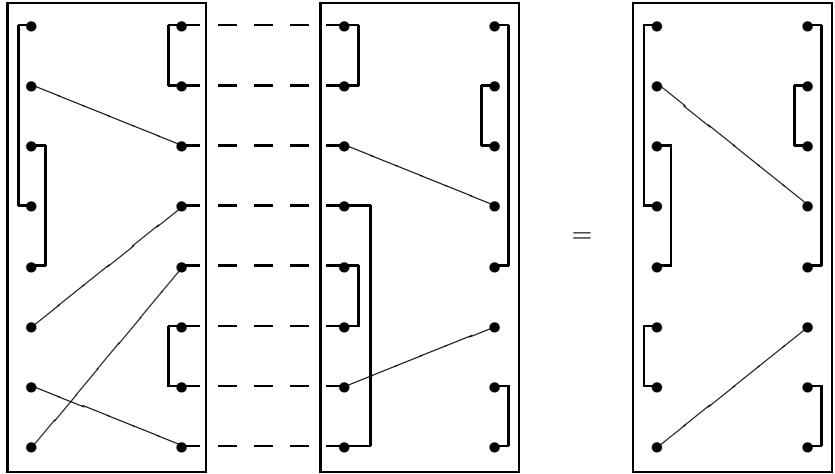
\begin{figure}
%%%%%%%%%%%%%%%%%%%%%%%%%%%%%%%%%%%%%%%%%%%%%%%%%%%%%%%%%%%%%%%%%%%%%%%%%%%%%%%%
%%%%%
%\input{figa1.pic}
%TexCad Options
%\grade{\on}
%\paths{\on}
%\beziermacro{\off}
%\reduce{\on}
%\snapping{\off}
%\quality{2.00}
%\graddiff{0.01}
%\snapasp{1}
%\zoom{1.00}
\special{em:linewidth 0.4pt} \unitlength 0.80mm
\linethickness{1pt}
\begin{picture}(150.00,75.00)
\put(13.00,00.00){\makebox(0,0)[cc]{$\bullet$}}
\put(13.00,10.00){\makebox(0,0)[cc]{$\bullet$}}
\put(13.00,20.00){\makebox(0,0)[cc]{$\bullet$}}
\put(13.00,30.00){\makebox(0,0)[cc]{$\bullet$}}
\put(13.00,40.00){\makebox(0,0)[cc]{$\bullet$}}
\put(13.00,50.00){\makebox(0,0)[cc]{$\bullet$}}
\put(13.00,60.00){\makebox(0,0)[cc]{$\bullet$}}
\put(13.00,70.00){\makebox(0,0)[cc]{$\bullet$}}
\put(38.00,10.00){\makebox(0,0)[cc]{$\bullet$}}
\put(38.00,20.00){\makebox(0,0)[cc]{$\bullet$}}
\put(38.00,30.00){\makebox(0,0)[cc]{$\bullet$}}
\put(38.00,40.00){\makebox(0,0)[cc]{$\bullet$}}
\put(38.00,50.00){\makebox(0,0)[cc]{$\bullet$}}
\put(38.00,60.00){\makebox(0,0)[cc]{$\bullet$}}
\put(38.00,70.00){\makebox(0,0)[cc]{$\bullet$}}
\put(38.00,00.00){\makebox(0,0)[cc]{$\bullet$}}

\put(65.00,00.00){\makebox(0,0)[cc]{$\bullet$}}
\put(65.00,10.00){\makebox(0,0)[cc]{$\bullet$}}
\put(65.00,20.00){\makebox(0,0)[cc]{$\bullet$}}
\put(65.00,30.00){\makebox(0,0)[cc]{$\bullet$}}
\put(65.00,40.00){\makebox(0,0)[cc]{$\bullet$}}
\put(65.00,50.00){\makebox(0,0)[cc]{$\bullet$}}
\put(65.00,60.00){\makebox(0,0)[cc]{$\bullet$}}
\put(65.00,70.00){\makebox(0,0)[cc]{$\bullet$}}
\put(90.00,10.00){\makebox(0,0)[cc]{$\bullet$}}
\put(90.00,20.00){\makebox(0,0)[cc]{$\bullet$}}
\put(90.00,30.00){\makebox(0,0)[cc]{$\bullet$}}
\put(90.00,40.00){\makebox(0,0)[cc]{$\bullet$}}
\put(90.00,50.00){\makebox(0,0)[cc]{$\bullet$}}
\put(90.00,60.00){\makebox(0,0)[cc]{$\bullet$}}
\put(90.00,70.00){\makebox(0,0)[cc]{$\bullet$}}
\put(90.00,00.00){\makebox(0,0)[cc]{$\bullet$}}
\drawline(38.00,00.30)(41.00,00.30)
\drawline(44.00,00.30)(47.00,00.30)
\drawline(50.00,00.30)(53.00,00.30)
\drawline(56.00,00.30)(59.00,00.30)
\drawline(62.00,00.30)(65.00,00.30)
\drawline(38.00,10.30)(41.00,10.30)
\drawline(44.00,10.30)(47.00,10.30)
\drawline(50.00,10.30)(53.00,10.30)
\drawline(56.00,10.30)(59.00,10.30)
\drawline(62.00,10.30)(65.00,10.30)
\drawline(38.00,20.30)(41.00,20.30)
\drawline(44.00,20.30)(47.00,20.30)
\drawline(50.00,20.30)(53.00,20.30)
\drawline(56.00,20.30)(59.00,20.30)
\drawline(62.00,20.30)(65.00,20.30)
\drawline(38.00,30.30)(41.00,30.30)
\drawline(44.00,30.30)(47.00,30.30)
\drawline(50.00,30.30)(53.00,30.30)
\drawline(56.00,30.30)(59.00,30.30)
\drawline(62.00,30.30)(65.00,30.30)
\drawline(38.00,40.30)(41.00,40.30)
\drawline(44.00,40.30)(47.00,40.30)
\drawline(50.00,40.30)(53.00,40.30)
\drawline(56.00,40.30)(59.00,40.30)
\drawline(62.00,40.30)(65.00,40.30)
\drawline(38.00,50.30)(41.00,50.30)
\drawline(44.00,50.30)(47.00,50.30)
\drawline(50.00,50.30)(53.00,50.30)
\drawline(56.00,50.30)(59.00,50.30)
\drawline(62.00,50.30)(65.00,50.30)
\drawline(38.00,60.30)(41.00,60.30)
\drawline(44.00,60.30)(47.00,60.30)
\drawline(50.00,60.30)(53.00,60.30)
\drawline(56.00,60.30)(59.00,60.30)
\drawline(62.00,60.30)(65.00,60.30)
\drawline(38.00,70.30)(41.00,70.30)
\drawline(44.00,70.30)(47.00,70.30)
\drawline(50.00,70.30)(53.00,70.30)
\drawline(56.00,70.30)(59.00,70.30)
\drawline(62.00,70.30)(65.00,70.30)
%%%%%%%%%%%%%%%%%%%%%%%%%%%%%%%%%%%%%%%%%%%%%%%%%%%%%%%%% FRAME %%%%%%%%%%%%%%%%%%%%%%%%%%%%%%%%

\drawline(09.00,-04.00)(09.00,74.00)
\drawline(09.00,74.00)(42.00,74.00)
\drawline(42.00,74.00)(42.00,-04.00)
\drawline(42.00,-04.00)(09.00,-04.00)

\drawline(61.00,-04.00)(61.00,74.00)
\drawline(61.00,74.00)(94.00,74.00)
\drawline(94.00,74.00)(94.00,-04.00)
\drawline(94.00,-04.00)(61.00,-04.00)

\drawline(113.00,-04.00)(113.00,74.00)
\drawline(113.00,74.00)(146.00,74.00)
\drawline(146.00,74.00)(146.00,-04.00)
\drawline(146.00,-04.00)(113.00,-04.00)

%%%%%%%%%%%%%%%%%%%%%%%%%%%%%%%%%%%%%%%%%%%%%%%%%%%%%%%%% FIRST ELEMENT %%%%%%%%%%%%%%%%%%%%%%%%%%%

\drawline(13.30,00.30)(38.30,30.30)
\drawline(13.30,10.30)(38.30,00.30)
\drawline(13.30,20.30)(38.30,40.30)
\drawline(13.30,60.30)(38.30,50.30)

\drawline(13.30,30.30)(15.30,30.30)
\drawline(15.30,30.30)(15.30,50.30)
\drawline(15.30,50.30)(13.30,50.30)

\drawline(13.30,40.30)(10.80,40.30)
\drawline(10.80,40.30)(10.80,70.30)
\drawline(10.80,70.30)(13.30,70.30)

\drawline(38.30,10.30)(35.80,10.30)
\drawline(35.80,10.30)(35.80,20.30)
\drawline(35.80,20.30)(38.30,20.30)

\drawline(38.30,60.30)(35.80,60.30)
\drawline(35.80,60.30)(35.80,70.30)
\drawline(35.80,70.30)(38.30,70.30)

%%%%%%%%%%%%%%%%%%%%%%%%%%%%%%%%%%%%%%%%%%%%%%%%%%%%%%%% SECOND ELEMENT %%%%%%%%%%%%%%%%%%%%%%%%%%%%

\drawline(65.30,10.30)(90.30,20.30)
\drawline(65.30,50.30)(90.30,40.30)

\drawline(65.30,00.30)(69.30,00.30)
\drawline(69.30,00.30)(69.30,40.30)
\drawline(69.30,40.30)(65.30,40.30)

\drawline(65.30,20.30)(67.30,20.30)
\drawline(67.30,20.30)(67.30,30.30)
\drawline(67.30,30.30)(65.30,30.30)

\drawline(65.30,70.30)(67.30,70.30)
\drawline(67.30,70.30)(67.30,60.30)
\drawline(67.30,60.30)(65.30,60.30)

\drawline(90.30,50.30)(87.80,50.30)
\drawline(87.80,50.30)(87.80,60.30)
\drawline(87.80,60.30)(90.30,60.30)

\drawline(90.30,00.30)(92.30,00.30)
\drawline(92.30,00.30)(92.30,10.30)
\drawline(92.30,10.30)(90.30,10.30)

\drawline(90.30,30.30)(92.30,30.30)
\drawline(92.30,30.30)(92.30,70.30)
\drawline(92.30,70.30)(90.30,70.30)

%%%%%%%%%%%%%%%%%%%%%%%%%%%%%%%%%%%%%%%%%%%%%% THE PRODUCT ELEMENT %%%%%%%%%%%%%%%%%%%%%%%%%%%%%

\drawline(117.30,00.30)(142.30,20.30)
\drawline(117.30,60.30)(142.30,40.30)

\drawline(117.30,10.30)(114.80,10.30)
\drawline(114.80,10.30)(114.80,20.30)
\drawline(114.80,20.30)(117.30,20.30)

\drawline(117.30,40.30)(114.80,40.30)
\drawline(114.80,40.30)(114.80,70.30)
\drawline(114.80,70.30)(117.30,70.30)

\drawline(117.30,30.30)(119.30,30.30)
\drawline(119.30,30.30)(119.30,50.30)
\drawline(119.30,50.30)(117.30,50.30)

\drawline(142.30,50.30)(139.80,50.30)
\drawline(139.80,50.30)(139.80,60.30)
\drawline(139.80,60.30)(142.30,60.30)

\drawline(142.30,00.30)(144.30,00.30)
\drawline(144.30,00.30)(144.30,10.30)
\drawline(144.30,10.30)(142.30,10.30)

\drawline(142.30,30.30)(144.30,30.30)
\drawline(144.30,30.30)(144.30,70.30)
\drawline(144.30,70.30)(142.30,70.30)

%%%%%%%%%%%%%%%%%%%%%%%%%%%%%%%%%%%%%%%%%%%%%%%%%%%%%%%%%%%%%%%%%%%%%%%%%%%%%%%%%%%%%%%%%
\put(117.00,00.00){\makebox(0,0)[cc]{$\bullet$}}
\put(117.00,10.00){\makebox(0,0)[cc]{$\bullet$}}
\put(117.00,20.00){\makebox(0,0)[cc]{$\bullet$}}
\put(117.00,30.00){\makebox(0,0)[cc]{$\bullet$}}
\put(117.00,40.00){\makebox(0,0)[cc]{$\bullet$}}
\put(117.00,50.00){\makebox(0,0)[cc]{$\bullet$}}
\put(117.00,60.00){\makebox(0,0)[cc]{$\bullet$}}
\put(117.00,70.00){\makebox(0,0)[cc]{$\bullet$}}
\put(142.00,10.00){\makebox(0,0)[cc]{$\bullet$}}
\put(142.00,20.00){\makebox(0,0)[cc]{$\bullet$}}
\put(142.00,30.00){\makebox(0,0)[cc]{$\bullet$}}
\put(142.00,40.00){\makebox(0,0)[cc]{$\bullet$}}
\put(142.00,50.00){\makebox(0,0)[cc]{$\bullet$}}
\put(142.00,60.00){\makebox(0,0)[cc]{$\bullet$}}
\put(142.00,70.00){\makebox(0,0)[cc]{$\bullet$}}
\put(142.00,00.00){\makebox(0,0)[cc]{$\bullet$}}
\put(104.50,35.00){\makebox(0,0)[cc]{$=$}}
\end{picture}
\caption{Elements of $\B_8$ and their
multiplication.}\label{fig:multiplication-Brauer}

%%%%%%%%%%%%%%%%%%%%%%%%%%%%%%%%%%%%%%%%%%%%%%%%%%%%%%%%%%%%%%%%%%%%%%%%%%%%%%%%
%%%%
\end{figure}

Now we would like to define the multiplication in $\B_n$. 
To give a formal definition, for $\pi\in\B_n$ and $x,y\in\bn\cup\bn'$ 
we set $x\equiv_{\pi}y$ provided that $x$ and $y$ are in the same block 
of $\pi$. The relation $\equiv_{\pi}$ is an equivalence relation on
$\bn\cup\bn'$ with two-element equivalence classes.
Take now $\pi,\tau\in\B_n$. Define a new equivalence relation, 
$\equiv$, on $\bn\cup\bn'$ as follows:
\begin{itemize}
\item
for $x,y\in\bn$ we have $x\equiv y$ if and only if $x\equiv_{\pi} y$ or
there is a sequence, $c_1,\dots,c_{2s}$, $s\geq 1$, of elements in
$\bn$, such that $x\equiv_{\pi} c'_1$, $c_1\equiv_{\tau} c_2$,
$c'_2\equiv_{\pi} c'_3$, \dots, $c_{2s-1}\equiv_{\tau} c_{2s}$, and
$c'_{2s}\equiv_{\pi} y$;
\item
for $x,y\in\bn$ we have $x'\equiv y'$ if and only if $x'\equiv_{\tau}
y'$ or there is a sequence, $c_1,\dots,c_{2s}$, $s\geq 1$, of
elements in $\bn$, such that  $x'\equiv_{\tau} c_1$, $c'_1\equiv_{\pi}
c'_2$, $c_2\equiv_{\tau} c_3$, \dots, $c'_{2s-1}\equiv_{\pi} c'_{2s}$,
and $c_{2s}\equiv_{\tau} y'$;
\item
for $x,y\in\bn$ we have $x\equiv y'$ if and only if $y'\equiv x$ if
and only if there is a sequence, $c_1,\dots$, $c_{2s-1}$, $s\geq
1$, of elements in $\bn$, such that $x\equiv_{\pi} c'_1$, $c_1\equiv_{\tau}
c_2$, $c'_2\equiv_{\pi} c'_3,\dots,$ $c'_{2s-2}\equiv_{\pi} c'_{2s-1}$,
and $c_{2s-1}\equiv_{\tau} y'$.
\end{itemize}
It is easy to see that $\equiv$ determines an equivalence relation
on $\bn\cup\bn'$ with two-element classes and thus is an element of 
$\B_n$. We define this element to be the product $\pi\tau$. It
is straightforward that this multiplication is associative.
In our geometric realization the above multiplication reduces to
concatenation of chips, see an example on
Figure~\ref{fig:multiplication-Brauer}.

Note that the element $\bigl\{\{k,k'\}_{k\in\bn}\bigr\}$ is the 
identity element in $\B_n$. It is easy to see (see for example 
\cite{Maz1}) that the group of all invertible elements in $\B_n$ is 
precisely the set of all elements of corank $0$, and it is isomorphic 
to $\Sym_n$. We identify the elements of this subgroup of $\B_n$ with 
$\Sym_n$ in the following way: $\pi\in \Sym_n$ corresponds to
the element $\bigl\{\{k,\pi(k)\}_{k\in\bn}\bigr\}$. Then the subsemigroup of 
all non-invertible elements of  $\B_n$ coincides with $\B_n\setminus\Sym_n$.

We denote by $\GR$, $\GL$, $\GH$, $\GD$ and $\GJ$ Green's relations, in 
particular, for a semigroup, $S$, and $a\in S$, $\mathcal{H}_a$ denotes
the $\mathcal{H}$-class of $S$  containing $a$ (similarly for all
other relations). We will need the following description of Green's 
relations for $\B_n$, which was obtained in \cite{Maz1}:

\begin{lemma}\label{lm:Green-on-Brauer}
Let $\pi,\tau\in\B_n$. Then
\begin{enumerate}[(i)]
\item $\pi\GR \tau$ if and only if $\pi$ and $\tau$ have the same 
left brackets;
\item $\pi\GL \tau$ if and only if $\pi$ and $\tau$ have the same 
right brackets;
\item $\pi\GH \tau$ if and only if $\pi$ and $\tau$ have both, the 
same left brackets and the same right brackets;
\item $\pi\GD \tau$ if and only if $\pi\GJ \tau$ if and only if $\corank(\pi)=\corank(\tau)$.
\end{enumerate}
\end{lemma}

\section{An irreducible system of generators for $\B_n\setminus\Sym_n$}\label{sec:irreducible-system}

\begin{figure}
%%%%%%%%%%%%%%%%%%%%%%%%%%%%%%%%%%%%%%%%%%%%%%%%%%%%%%%%%%%%%%%%%%%%%%%%%%%%%%%%
%%%%%
%\input{figa1.pic}
%TexCad Options
%\grade{\on}
%\paths{\on}
%\beziermacro{\off}
%\reduce{\on}
%\snapping{\off}
%\quality{2.00}
%\graddiff{0.01}
%\snapasp{1}
%\zoom{1.00}
\special{em:linewidth 0.4pt} \unitlength 0.80mm
\linethickness{1pt}
\begin{picture}(150.00,35.00)
\put(71.00,00.00){\makebox(0,0)[cc]{$\bullet$}}
\put(71.00,10.00){\makebox(0,0)[cc]{$\bullet$}}
\put(71.00,20.00){\makebox(0,0)[cc]{$\bullet$}}
\put(71.00,30.00){\makebox(0,0)[cc]{$\bullet$}}
\put(81.00,00.00){\makebox(0,0)[cc]{$\bullet$}}
\put(81.00,10.00){\makebox(0,0)[cc]{$\bullet$}}
\put(81.00,20.00){\makebox(0,0)[cc]{$\bullet$}}
\put(81.00,30.00){\makebox(0,0)[cc]{$\bullet$}}

\drawline(67.00,-04.00)(67.00,34.00)
\drawline(67.00,34.00)(85.00,34.00)
\drawline(85.00,34.00)(85.00,-04.00)
\drawline(85.00,-04.00)(67.00,-04.00)

\drawline(71.30,20.30)(81.30,20.30)
\drawline(71.30,00.30)(81.30,00.30)
                     
\drawline(71.30,30.30)(73.30,30.30)
\drawline(73.30,30.30)(73.30,10.30)
\drawline(73.30,10.30)(71.30,10.30)                  

\drawline(81.30,30.30)(78.80,30.30)
\drawline(78.80,30.30)(78.80,10.30)
\drawline(78.80,10.30)(81.30,10.30)
\end{picture}
\caption{The atom $\sigma_{1,3}$ of $\B_4$.}
\label{fig:example-for-atom}
%%%%%%%%%%%%%%%%%%%%%%%%%%%%%%%%%%%%%%%%%%%%%%%%%%%%%%%%%%%%%%
\end{figure}
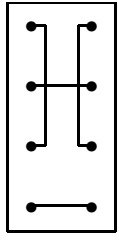

For $i,j\in\bn$, $i\neq j$, define $\sigma_{i,j}$ as follows:
\begin{equation*}
\sigma_{i,j}=\bigl\{\{i,j\},\{i',j'\},\{k,k'\}_{k\ne i,j}\bigr\}.
\end{equation*}
We have $\sigma_{i,j}=\sigma_{j,i}=\sigma_{i,j}^2$ and
$\corank(\sigma_{i,j})=2$. We will call these elements \emph{atoms}. 
An example of an atom can be found on
Figure~\ref{fig:example-for-atom}.

The following statement was proved in \cite{Mal}. However, because of
the poor availability of \cite{Mal} we will prove it here as well.

\begin{proposition}\label{pr:irreducible-system}
The set of all atoms is an irreducible system of generators 
in $\B_n\setminus\Sym_n$.
\end{proposition}

To prove this statement we will need several auxiliary lemmas.

\begin{lemma}\label{lm:bn-sn-generated-by-cr-2}
The semigroup $\B_n\setminus\Sym_n$ is 
generated by the set of all elements of 
corank $2$.
\end{lemma}

\begin{proof}
Let $\pi\in \B_n\setminus\Sym_n$ be written in the form 
\eqref{eq:convenient-form} as follows:
\begin{equation*}
\pi=\bigl\{\{i,\theta(i)'\}_{i\in I},\{u_j,v_j\}_{j\in J},\{f_j',g_j'\}_{j\in J}\bigr\}.
\end{equation*}
We have $J\ne\varnothing$, $I\subset \bn$,  and  $\theta: I\to \bn$ 
is an injection. Assume that $\corank(\pi)>2$. Fix $j_0\in J$. 
Construct a bijection, $\vartheta:{\bn}\setminus\{u_{j_0},v_{j_0}\}
\to{\bn}\setminus\{f_{j_0},g_{j_0}\}$ as follows:
\begin{itemize}
\item
$\vartheta(i)=\theta(i)$ for all $i\in I$;
\item
$\vartheta(u_j)=f_j$ and $\vartheta(v_j)=g_j$ for all $j\in J\setminus\{j_0\}$.
\end{itemize}
Put now 
$\tau=
\bigl\{\{i,\vartheta(i)'\}_{i\neq u_{j_0},v_{j_0}},\{u_{j_0},v_{j_0}\},
\{f_{j_0}',g_{j_0}'\}\bigr\}$.
We have $\corank(\tau)=2$ and a direct calculation shows 
that $\pi=\prod\limits_{j\in J} \sigma_{u_j,v_j}\cdot \tau$. 
The statement follows.
\end{proof}

\begin{lemma}\label{lm:gr-el-is-decomposable}
Every element of the maximal subgroup, corresponding to an atom, 
is decomposable into a product of atoms.
\end{lemma}

\begin{proof}
Let $\pi\in\B_n$ be a group element of corank $2$, $\GH$-related to some atom. From 
Lemma~\ref{lm:Green-on-Brauer} it follows that in this case
$\pi=\bigl\{\{i,\theta(i)'\}_{i\ne u,v},\{u,v\},\{u',v'\}\bigr\}$ 
for some $u,v\in{\bn}$, $u\neq v$, and some bijection, 
$\theta:{\bn}\setminus\{u,v\}\to {\bn}\setminus\{u,v\}$. 
We consider $\theta$ as an element of $\Sym_{{\bn}\setminus\{u,v\}}$. 
Let 
\begin{displaymath}
\theta=(i_1^{(1)},\dots,i_{p_1}^{(1)})\cdot\cdots\cdot
(i_1^{(s)},\dots,i_{p_s}^{(s)})
\end{displaymath}
be a cyclic decomposition of $\theta$.
By a direct calculation one obtains that
\begin{equation}\label{cycleformula}
\pi=\sigma_{u,v}\sigma_{u,i_{1}^{(1)}}\dots
\sigma_{u,i_{p_1}^{(1)}}\sigma_{u,v}\cdot\dots\cdot
\sigma_{u,v}\sigma_{u,i_{1}^{(s)}}\dots
\sigma_{u,i_{p_s}^{(s)}}\sigma_{u,v}.
\end{equation}
The statement follows.
\end{proof}

Now we are ready to prove Proposition~\ref{pr:irreducible-system}:

\begin{proof}[Proof of Proposition~\ref{pr:irreducible-system}.]
First we show that atoms generate $\B_n\setminus\Sym_n$. Because of
Lemma~\ref{lm:bn-sn-generated-by-cr-2} it is enough to show that
any element $\pi\in\B_n$ of corank $2$ decomposes into a product
of atoms. We again write $\pi$ in the form \eqref{eq:convenient-form}:
\begin{displaymath}
\pi=\bigl\{\{i,\theta (i)'\}_{i\in I},\{u,v\},\{f',g'\}\bigr\},
\end{displaymath}
where $u,v\in\bn$, $u\neq v$; $f,g\in \bn$, $f\neq g$;
and $\theta:\bn\setminus \{u,v\}\to\bn\setminus \{f,g\}$ is a bijection.
Without loss of generality we may assume that $v\ne f$. 
Consider the element $\tau=\sigma_{v,f}\sigma_{f,g}=
\bigl\{\{v,f\},\{f',g'\},\{g,v'\},\{k,k'\}_{k\ne v,f,g}\bigr\}$. 
From Lemma~\ref{lm:Green-on-Brauer} we have $\pi\GH\sigma_{u,v}\tau$ and
$\sigma_{u,v}\GR\sigma_{u,v}\tau$.
Hence, due to Green's Lemma, 
we have  that the map $x\mapsto x\tau$ from $\mathcal{H}_{\sigma_{u,v}}$ 
to $\mathcal{H}_{\sigma_{u,v}\tau}$  is a bijection. Therefore there exists 
$\xi\in \mathcal{H}_{\sigma_{u,v}}$ such 
that $\pi=\xi\tau$. By Lemma~\ref{lm:gr-el-is-decomposable}, $\xi$  
decomposes into a product of atoms. Hence so does $\pi$ as well.

Now we prove that no atom can be decomposed into a product of other
atoms. Let $\sigma_{u,v}=\sigma_{u_1,v_1}\dots \sigma_{u_k,v_k}$.
The product $\sigma_{u_1,v_1}\dots \sigma_{u_k,v_k}$ must contain
the left bracket $\{u_1,v_1\}$ by the definition of the 
multiplication in $\B_n$. However, the element $\sigma_{u,v}$ contains
the unique left bracket $\{u,v\}$. This implies that 
$\{u_1,v_1\}=\{u,v\}$ and the necessary statement follows. The proof is 
complete.
\end{proof}

After Proposition~\ref{pr:irreducible-system} it is natural to ask
what is the presentation of $\B_n\setminus\Sym_n$ with respect 
to the system $\{\sigma_{u,v}\}$ of generators. We answer this
question in the next section.

%%%%%%%%%%%%%%%%%%%%%%%%%%%%%%%%%%%%%%%%%%%%%%%%%%%%%%%%%%%%%%%%%%%%%%%%%%%%%%%%%%%%%%%%%%%%

\section{Main result}\label{sec:main}

Denote by $T$ the semigroup generated by $\tau_{i,j}$, $i,j\in\bn$,
$i\neq j$; subject to the following relations (here $i$, $j$, $k$, $l$
are pairwise different):
\begin{eqnarray}
\tau_{i,j}&=&\tau_{j,i};\label{1}\\
\tau_{i,j}^2&=&\tau_{i,j};\label{2}\\
\tau_{i,j}\tau_{j,k}\tau_{k,l}&=&\tau_{i,j}\tau_{i,l}\tau_{k,l};\label{3}\\
\tau_{i,j}\tau_{i,k}\tau_{j,k}&=&\tau_{i,j}\tau_{j,k};\label{4}\\
\tau_{i,j}\tau_{j,k}\tau_{i,j}&=&\tau_{i,j};\label{5}\\
\tau_{i,j}\tau_{k,l}\tau_{i,k}&=&\tau_{i,j}\tau_{j,l}\tau_{i,k};\label{6}\\
\tau_{i,j}\tau_{k,l}&=&\tau_{k,l}\tau_{i,j}\label{7}.
\end{eqnarray}

\begin{figure}
%%%%%%%%%%%%%%%%%%%%%%%%%%%%%%%%%%%%%%%%%%%%%%%%%%%%%%%%%%%%%%%%%%%%%%%%%%%%%%%%
%%%%%
%\input{figa1.pic}
%TexCad Options
%\grade{\on}
%\paths{\on}
%\beziermacro{\off}
%\reduce{\on}
%\snapping{\off}
%\quality{2.00}
%\graddiff{0.01}
%\snapasp{1}
%\zoom{1.00}
\special{em:linewidth 0.4pt} \unitlength 0.80mm
\linethickness{1pt}
\begin{picture}(150.00,25.00)
%%%%%%%%%%%%%%%%%%%%%%%%%%%%%%%%%%%%%%%%%%%%%%%% FIRST LEFT ELEMENT %%%%%%%%%%%%%%%%%%%%%%%%%%%%%%%%%%%
\put(11.00,00.00){\makebox(0,0)[cc]{$\bullet$}}
\put(11.00,10.00){\makebox(0,0)[cc]{$\bullet$}}
\put(11.00,20.00){\makebox(0,0)[cc]{$\bullet$}}
\put(21.00,00.00){\makebox(0,0)[cc]{$\bullet$}}
\put(21.00,10.00){\makebox(0,0)[cc]{$\bullet$}}
\put(21.00,20.00){\makebox(0,0)[cc]{$\bullet$}}

\drawline(07.00,-04.00)(07.00,24.00)
\drawline(07.00,24.00)(25.00,24.00)
\drawline(25.00,24.00)(25.00,-04.00)
\drawline(25.00,-04.00)(07.00,-04.00)

\drawline(11.30,00.30)(21.30,00.30)

\drawline(11.30,20.30)(13.30,20.30)
\drawline(13.30,20.30)(13.30,10.30)
\drawline(13.30,10.30)(11.30,10.30)

\drawline(21.30,20.30)(18.80,20.30)
\drawline(18.80,20.30)(18.80,10.30)
\drawline(18.80,10.30)(21.30,10.30)

%%%%%%%%%%%%%%%%%%%%%%%%%%%%%%%%%%%%%%%%%%%%%%%%%%% SECOND LEFT ELEMENT %%%%%%%%%%%%%%%%%%%%%%%%%%%%

\put(34.00,00.00){\makebox(0,0)[cc]{$\bullet$}}
\put(34.00,10.00){\makebox(0,0)[cc]{$\bullet$}}
\put(34.00,20.00){\makebox(0,0)[cc]{$\bullet$}}
\put(44.00,00.00){\makebox(0,0)[cc]{$\bullet$}}
\put(44.00,10.00){\makebox(0,0)[cc]{$\bullet$}}
\put(44.00,20.00){\makebox(0,0)[cc]{$\bullet$}}

\drawline(30.00,-04.00)(30.00,24.00)
\drawline(30.00,24.00)(48.00,24.00)
\drawline(48.00,24.00)(48.00,-04.00)
\drawline(48.00,-04.00)(30.00,-04.00)

\drawline(34.30,10.30)(44.30,10.30)

\drawline(34.30,00.30)(36.30,00.30)
\drawline(36.30,00.30)(36.30,20.30)
\drawline(36.30,20.30)(34.30,20.30)

\drawline(44.30,00.30)(41.80,00.30)
\drawline(41.80,00.30)(41.80,20.30)
\drawline(41.80,20.30)(44.30,20.30)

%%%%%%%%%%%%%%%%%%%%%%%%%%%%%%%%%%%%%%%%%%%%%%%%%%% THIRD LEFT ELEMENT %%%%%%%%%%%%%%%%%%%%%%%%%%%%

\put(57.00,00.00){\makebox(0,0)[cc]{$\bullet$}}
\put(57.00,10.00){\makebox(0,0)[cc]{$\bullet$}}
\put(57.00,20.00){\makebox(0,0)[cc]{$\bullet$}}
\put(67.00,00.00){\makebox(0,0)[cc]{$\bullet$}}
\put(67.00,10.00){\makebox(0,0)[cc]{$\bullet$}}
\put(67.00,20.00){\makebox(0,0)[cc]{$\bullet$}}

\drawline(53.00,-04.00)(53.00,24.00)
\drawline(53.00,24.00)(71.00,24.00)
\drawline(71.00,24.00)(71.00,-04.00)
\drawline(71.00,-04.00)(53.00,-04.00)

\drawline(57.30,20.30)(67.30,20.30)

\drawline(57.30,00.30)(59.30,00.30)
\drawline(59.30,00.30)(59.30,10.30)
\drawline(59.30,10.30)(57.30,10.30)

\drawline(67.30,00.30)(64.80,00.30)
\drawline(64.80,00.30)(64.80,10.30)
\drawline(64.80,10.30)(67.30,10.30)

%%%%%%%%%%%%%%%%%%%%%%%%%%%%%%%%%%%%%%%%%%%%%%%%%%% FIRST RIGHT ELEMENT %%%%%%%%%%%%%%%%%%%%%%%%%%%%

\put(85.00,00.00){\makebox(0,0)[cc]{$\bullet$}}
\put(85.00,10.00){\makebox(0,0)[cc]{$\bullet$}}
\put(85.00,20.00){\makebox(0,0)[cc]{$\bullet$}}
\put(95.00,00.00){\makebox(0,0)[cc]{$\bullet$}}
\put(95.00,10.00){\makebox(0,0)[cc]{$\bullet$}}
\put(95.00,20.00){\makebox(0,0)[cc]{$\bullet$}}

\drawline(81.00,-04.00)(81.00,24.00)
\drawline(81.00,24.00)(99.00,24.00)
\drawline(99.00,24.00)(99.00,-04.00)
\drawline(99.00,-04.00)(81.00,-04.00)

\drawline(85.30,00.30)(95.30,00.30)

\drawline(85.30,20.30)(87.30,20.30)
\drawline(87.30,20.30)(87.30,10.30)
\drawline(87.30,10.30)(85.30,10.30)

\drawline(95.30,20.30)(92.80,20.30)
\drawline(92.80,20.30)(92.80,10.30)
\drawline(92.80,10.30)(95.30,10.30)

%%%%%%%%%%%%%%%%%%%%%%%%%%%%%%%%%%%%%%%%%%%%%%%%%%% SECOND RIGHT ELEMENT %%%%%%%%%%%%%%%%%%%%%%%%%%%%

\put(108.00,00.00){\makebox(0,0)[cc]{$\bullet$}}
\put(108.00,10.00){\makebox(0,0)[cc]{$\bullet$}}
\put(108.00,20.00){\makebox(0,0)[cc]{$\bullet$}}
\put(118.00,00.00){\makebox(0,0)[cc]{$\bullet$}}
\put(118.00,10.00){\makebox(0,0)[cc]{$\bullet$}}
\put(118.00,20.00){\makebox(0,0)[cc]{$\bullet$}}

\drawline(104.00,-04.00)(104.00,24.00)
\drawline(104.00,24.00)(122.00,24.00)
\drawline(122.00,24.00)(122.00,-04.00)
\drawline(122.00,-04.00)(104.00,-04.00)

\drawline(108.30,20.30)(118.30,20.30)

\drawline(108.30,00.30)(110.30,00.30)
\drawline(110.30,00.30)(110.30,10.30)
\drawline(110.30,10.30)(108.30,10.30)

\drawline(118.30,00.30)(115.80,00.30)
\drawline(115.80,00.30)(115.80,10.30)
\drawline(115.80,10.30)(118.30,10.30)

%%%%%%%%%%%%%%%%%%%%%%%%%%%%%%%%%%%%%%%%%%%%%%%%%% 1-2 %%%%%%%%%%%%%%%%%%%%%%%%%%%%%%%%%%%%%%%%%%%%%%%

\drawline(21.30,00.30)(24.30,00.30)
\drawline(26.30,00.30)(29.30,00.30)
\drawline(31.30,00.30)(34.30,00.30)

\drawline(21.30,10.30)(24.30,10.30)
\drawline(26.30,10.30)(29.30,10.30)
\drawline(31.30,10.30)(34.30,10.30)

\drawline(21.30,20.30)(24.30,20.30)
\drawline(26.30,20.30)(29.30,20.30)
\drawline(31.30,20.30)(34.30,20.30)

%%%%%%%%%%%%%%%%%%%%%%%%%%%%%%%%%%%%%%%%%%%%%%%%%% 2-3 %%%%%%%%%%%%%%%%%%%%%%%%%%%%%%%%%%%%%%%%%%%%%%%

\drawline(44.30,00.30)(47.30,00.30)
\drawline(49.30,00.30)(52.30,00.30)
\drawline(54.30,00.30)(57.30,00.30)

\drawline(44.30,10.30)(47.30,10.30)
\drawline(49.30,10.30)(52.30,10.30)
\drawline(54.30,10.30)(57.30,10.30)

\drawline(44.30,20.30)(47.30,20.30)
\drawline(49.30,20.30)(52.30,20.30)
\drawline(54.30,20.30)(57.30,20.30)

%%%%%%%%%%%%%%%%%%%%%%%%%%%%%%%%%%%%%%%%%%%%%%%%%% 4-5 %%%%%%%%%%%%%%%%%%%%%%%%%%%%%%%%%%%%%%%%%%%%%%%

\drawline(95.30,00.30)(98.30,00.30)
\drawline(100.30,00.30)(103.30,00.30)
\drawline(105.30,00.30)(108.30,00.30)

\drawline(95.30,10.30)(98.30,10.30)
\drawline(100.30,10.30)(103.30,10.30)
\drawline(105.30,10.30)(108.30,10.30)

\drawline(95.30,20.30)(98.30,20.30)
\drawline(100.30,20.30)(103.30,20.30)
\drawline(105.30,20.30)(108.30,20.30)

%%%%%%%%%%%%%%%%%%%%%%%%%%%%%%%%%%%%%%%%%%%%%%%%%%%%%%%%%%%%%%%%%%%%%%%%%%%%%%%%%%%%%%%%%%%%%%%%%%

%%%%%%%%%%%%%%%%%%%%%%%%%%%%%%%%%%%%%%%%%%%%%%%%%%%%%%%%%%%%%%%%%%%%%%%%%%%%%%%%%%%%%%%%%
\put(76.00,10.00){\makebox(0,0)[cc]{$=$}}

\put(04.30,00.30){\makebox(0,0)[cc]{$k$}}
\put(04.30,10.30){\makebox(0,0)[cc]{$j$}}
\put(04.30,20.30){\makebox(0,0)[cc]{$i$}}
\end{picture}
\caption{An example illustrating the
  relation~\eqref{4}.}\label{fig:example-for-4}

%%%%%%%%%%%%%%%%%%%%%%%%%%%%%%%%%%%%%%%%%%%%%%%%%%%%%%%%%%%%%%%%%%%%%%%%%%%%%%%%
%%%%
\end{figure}
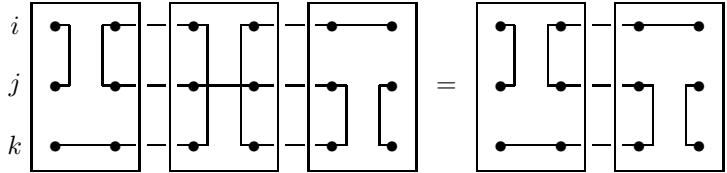

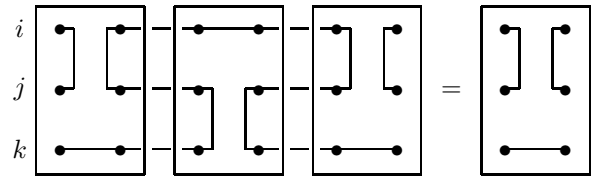
\begin{figure}
%%%%%%%%%%%%%%%%%%%%%%%%%%%%%%%%%%%%%%%%%%%%%%%%%%%%%%%%%%%%%%%%%%%%%%%%%%%%%%%%
%%%%%
%\input{figa1.pic}
%TexCad Options
%\grade{\on}
%\paths{\on}
%\beziermacro{\off}
%\reduce{\on}
%\snapping{\off}
%\quality{2.00}
%\graddiff{0.01}
%\snapasp{1}
%\zoom{1.00}
\special{em:linewidth 0.4pt} \unitlength 0.80mm
\linethickness{1pt}
\begin{center}
\begin{picture}(150.00,25.00)
%%%%%%%%%%%%%%%%%%%%%%%%%%%%%%%%%%%%%%%%%%%%%%%% FIRST LEFT ELEMENT %%%%%%%%%%%%%%%%%%%%%%%%%%%%%%%%%%%
\put(11.00,00.00){\makebox(0,0)[cc]{$\bullet$}}
\put(11.00,10.00){\makebox(0,0)[cc]{$\bullet$}}
\put(11.00,20.00){\makebox(0,0)[cc]{$\bullet$}}
\put(21.00,00.00){\makebox(0,0)[cc]{$\bullet$}}
\put(21.00,10.00){\makebox(0,0)[cc]{$\bullet$}}
\put(21.00,20.00){\makebox(0,0)[cc]{$\bullet$}}

\drawline(07.00,-04.00)(07.00,24.00)
\drawline(07.00,24.00)(25.00,24.00)
\drawline(25.00,24.00)(25.00,-04.00)
\drawline(25.00,-04.00)(07.00,-04.00)

\drawline(11.30,00.30)(21.30,00.30)

\drawline(11.30,20.30)(13.30,20.30)
\drawline(13.30,20.30)(13.30,10.30)
\drawline(13.30,10.30)(11.30,10.30)

\drawline(21.30,20.30)(18.80,20.30)
\drawline(18.80,20.30)(18.80,10.30)
\drawline(18.80,10.30)(21.30,10.30)

%%%%%%%%%%%%%%%%%%%%%%%%%%%%%%%%%%%%%%%%%%%%%%%%%%% SECOND LEFT ELEMENT %%%%%%%%%%%%%%%%%%%%%%%%%%%%

\put(34.00,00.00){\makebox(0,0)[cc]{$\bullet$}}
\put(34.00,10.00){\makebox(0,0)[cc]{$\bullet$}}
\put(34.00,20.00){\makebox(0,0)[cc]{$\bullet$}}
\put(44.00,00.00){\makebox(0,0)[cc]{$\bullet$}}
\put(44.00,10.00){\makebox(0,0)[cc]{$\bullet$}}
\put(44.00,20.00){\makebox(0,0)[cc]{$\bullet$}}

\drawline(30.00,-04.00)(30.00,24.00)
\drawline(30.00,24.00)(48.00,24.00)
\drawline(48.00,24.00)(48.00,-04.00)
\drawline(48.00,-04.00)(30.00,-04.00)

\drawline(34.30,20.30)(44.30,20.30)

\drawline(34.30,00.30)(36.30,00.30)
\drawline(36.30,00.30)(36.30,10.30)
\drawline(36.30,10.30)(34.30,10.30)

\drawline(44.30,00.30)(41.80,00.30)
\drawline(41.80,00.30)(41.80,10.30)
\drawline(41.80,10.30)(44.30,10.30)

%%%%%%%%%%%%%%%%%%%%%%%%%%%%%%%%%%%%%%%%%%%%%%%%%%% THIRD LEFT ELEMENT %%%%%%%%%%%%%%%%%%%%%%%%%%%%

\put(57.00,00.00){\makebox(0,0)[cc]{$\bullet$}}
\put(57.00,10.00){\makebox(0,0)[cc]{$\bullet$}}
\put(57.00,20.00){\makebox(0,0)[cc]{$\bullet$}}
\put(67.00,00.00){\makebox(0,0)[cc]{$\bullet$}}
\put(67.00,10.00){\makebox(0,0)[cc]{$\bullet$}}
\put(67.00,20.00){\makebox(0,0)[cc]{$\bullet$}}

\drawline(53.00,-04.00)(53.00,24.00)
\drawline(53.00,24.00)(71.00,24.00)
\drawline(71.00,24.00)(71.00,-04.00)
\drawline(71.00,-04.00)(53.00,-04.00)

\drawline(57.30,00.30)(67.30,00.30)

\drawline(57.30,20.30)(59.30,20.30)
\drawline(59.30,20.30)(59.30,10.30)
\drawline(59.30,10.30)(57.30,10.30)

\drawline(67.30,20.30)(64.80,20.30)
\drawline(64.80,20.30)(64.80,10.30)
\drawline(64.80,10.30)(67.30,10.30)

%%%%%%%%%%%%%%%%%%%%%%%%%%%%%%%%%%%%%%%%%%%%%%%%%%% FIRST RIGHT ELEMENT %%%%%%%%%%%%%%%%%%%%%%%%%%%%

\put(85.00,00.00){\makebox(0,0)[cc]{$\bullet$}}
\put(85.00,10.00){\makebox(0,0)[cc]{$\bullet$}}
\put(85.00,20.00){\makebox(0,0)[cc]{$\bullet$}}
\put(95.00,00.00){\makebox(0,0)[cc]{$\bullet$}}
\put(95.00,10.00){\makebox(0,0)[cc]{$\bullet$}}
\put(95.00,20.00){\makebox(0,0)[cc]{$\bullet$}}

\drawline(81.00,-04.00)(81.00,24.00)
\drawline(81.00,24.00)(99.00,24.00)
\drawline(99.00,24.00)(99.00,-04.00)
\drawline(99.00,-04.00)(81.00,-04.00)

\drawline(85.30,00.30)(95.30,00.30)

\drawline(85.30,20.30)(87.30,20.30)
\drawline(87.30,20.30)(87.30,10.30)
\drawline(87.30,10.30)(85.30,10.30)

\drawline(95.30,20.30)(92.80,20.30)
\drawline(92.80,20.30)(92.80,10.30)
\drawline(92.80,10.30)(95.30,10.30)

%%%%%%%%%%%%%%%%%%%%%%%%%%%%%%%%%%%%%%%%%%%%%%%%%% 1-2 %%%%%%%%%%%%%%%%%%%%%%%%%%%%%%%%%%%%%%%%%%%%%%%

\drawline(21.30,00.30)(24.30,00.30)
\drawline(26.30,00.30)(29.30,00.30)
\drawline(31.30,00.30)(34.30,00.30)

\drawline(21.30,10.30)(24.30,10.30)
\drawline(26.30,10.30)(29.30,10.30)
\drawline(31.30,10.30)(34.30,10.30)

\drawline(21.30,20.30)(24.30,20.30)
\drawline(26.30,20.30)(29.30,20.30)
\drawline(31.30,20.30)(34.30,20.30)

%%%%%%%%%%%%%%%%%%%%%%%%%%%%%%%%%%%%%%%%%%%%%%%%%% 2-3 %%%%%%%%%%%%%%%%%%%%%%%%%%%%%%%%%%%%%%%%%%%%%%%

\drawline(44.30,00.30)(47.30,00.30)
\drawline(49.30,00.30)(52.30,00.30)
\drawline(54.30,00.30)(57.30,00.30)

\drawline(44.30,10.30)(47.30,10.30)
\drawline(49.30,10.30)(52.30,10.30)
\drawline(54.30,10.30)(57.30,10.30)

\drawline(44.30,20.30)(47.30,20.30)
\drawline(49.30,20.30)(52.30,20.30)
\drawline(54.30,20.30)(57.30,20.30)

%%%%%%%%%%%%%%%%%%%%%%%%%%%%%%%%%%%%%%%%%%%%%%%%%%%%%%%%%%%%%%%%%%%%%%%%%%%%%%%%%%%%%%%%%
\put(76.00,10.00){\makebox(0,0)[cc]{$=$}}
\put(04.30,00.30){\makebox(0,0)[cc]{$k$}}
\put(04.30,10.30){\makebox(0,0)[cc]{$j$}}
\put(04.30,20.30){\makebox(0,0)[cc]{$i$}}
\end{picture}
\end{center}
\caption{An example illustrating the
 relation~\eqref{5}.}\label{fig:example-for-5}

%%%%%%%%%%%%%%%%%%%%%%%%%%%%%%%%%%%%%%%%%%%%%%%%%%%%%%%%%%%%%%%%%%%%%%%%%%%%%%%%
%%%%
\end{figure}

\begin{figure}
%%%%%%%%%%%%%%%%%%%%%%%%%%%%%%%%%%%%%%%%%%%%%%%%%%%%%%%%%%%%%%%%%%%%%%%%%%%%%%%%
%%%%%
%\input{figa1.pic}
%TexCad Options
%\grade{\on}
%\paths{\on}
%\beziermacro{\off}
%\reduce{\on}
%\snapping{\off}
%\quality{2.00}
%\graddiff{0.01}
%\snapasp{1}
%\zoom{1.00}
\special{em:linewidth 0.4pt} \unitlength 0.80mm
\linethickness{1pt}
\begin{picture}(150.00,35.00)
%%%%%%%%%%%%%%%%%%%%%%%%%%%%%%%%%%%%%%%%%%%%%%%% FIRST LEFT ELEMENT %%%%%%%%%%%%%%%%%%%%%%%%%%%%%%%%%%%
\put(11.00,00.00){\makebox(0,0)[cc]{$\bullet$}}
\put(11.00,10.00){\makebox(0,0)[cc]{$\bullet$}}
\put(11.00,20.00){\makebox(0,0)[cc]{$\bullet$}}
\put(11.00,30.00){\makebox(0,0)[cc]{$\bullet$}}
\put(21.00,00.00){\makebox(0,0)[cc]{$\bullet$}}
\put(21.00,10.00){\makebox(0,0)[cc]{$\bullet$}}
\put(21.00,20.00){\makebox(0,0)[cc]{$\bullet$}}
\put(21.00,30.00){\makebox(0,0)[cc]{$\bullet$}}

\drawline(07.00,-04.00)(07.00,34.00)
\drawline(07.00,34.00)(25.00,34.00)
\drawline(25.00,34.00)(25.00,-04.00)
\drawline(25.00,-04.00)(07.00,-04.00)

\drawline(11.30,00.30)(21.30,00.30)
\drawline(11.30,10.30)(21.30,10.30)

\drawline(11.30,20.30)(13.30,20.30)
\drawline(13.30,20.30)(13.30,30.30)
\drawline(13.30,30.30)(11.30,30.30)

\drawline(21.30,20.30)(18.80,20.30)
\drawline(18.80,20.30)(18.80,30.30)
\drawline(18.80,30.30)(21.30,30.30)

%%%%%%%%%%%%%%%%%%%%%%%%%%%%%%%%%%%%%%%%%%%%%%%%%%% SECOND LEFT ELEMENT %%%%%%%%%%%%%%%%%%%%%%%%%%%%

\put(34.00,00.00){\makebox(0,0)[cc]{$\bullet$}}
\put(34.00,10.00){\makebox(0,0)[cc]{$\bullet$}}
\put(34.00,20.00){\makebox(0,0)[cc]{$\bullet$}}
\put(34.00,30.00){\makebox(0,0)[cc]{$\bullet$}}
\put(44.00,00.00){\makebox(0,0)[cc]{$\bullet$}}
\put(44.00,10.00){\makebox(0,0)[cc]{$\bullet$}}
\put(44.00,20.00){\makebox(0,0)[cc]{$\bullet$}}
\put(44.00,30.00){\makebox(0,0)[cc]{$\bullet$}}

\drawline(30.00,-04.00)(30.00,34.00)
\drawline(30.00,34.00)(48.00,34.00)
\drawline(48.00,34.00)(48.00,-04.00)
\drawline(48.00,-04.00)(30.00,-04.00)

\drawline(34.30,20.30)(44.30,20.30)
\drawline(34.30,30.30)(44.30,30.30)

\drawline(34.30,10.30)(36.30,10.30)
\drawline(36.30,10.30)(36.30,00.30)
\drawline(36.30,00.30)(34.30,00.30)

\drawline(44.30,10.30)(41.80,10.30)
\drawline(41.80,10.30)(41.80,00.30)
\drawline(41.80,00.30)(44.30,00.30)

%%%%%%%%%%%%%%%%%%%%%%%%%%%%%%%%%%%%%%%%%%%%%%%%%%% THIRD LEFT ELEMENT %%%%%%%%%%%%%%%%%%%%%%%%%%%%

\put(57.00,00.00){\makebox(0,0)[cc]{$\bullet$}}
\put(57.00,10.00){\makebox(0,0)[cc]{$\bullet$}}
\put(57.00,20.00){\makebox(0,0)[cc]{$\bullet$}}
\put(57.00,30.00){\makebox(0,0)[cc]{$\bullet$}}
\put(67.00,00.00){\makebox(0,0)[cc]{$\bullet$}}
\put(67.00,10.00){\makebox(0,0)[cc]{$\bullet$}}
\put(67.00,20.00){\makebox(0,0)[cc]{$\bullet$}}
\put(67.00,30.00){\makebox(0,0)[cc]{$\bullet$}}

\drawline(53.00,-04.00)(53.00,34.00)
\drawline(53.00,34.00)(71.00,34.00)
\drawline(71.00,34.00)(71.00,-04.00)
\drawline(71.00,-04.00)(53.00,-04.00)

\drawline(57.30,20.30)(67.30,20.30)
\drawline(57.30,00.30)(67.30,00.30)

\drawline(57.30,30.30)(59.30,30.30)
\drawline(59.30,30.30)(59.30,10.30)
\drawline(59.30,10.30)(57.30,10.30)

\drawline(67.30,30.30)(64.80,30.30)
\drawline(64.80,30.30)(64.80,10.30)
\drawline(64.80,10.30)(67.30,10.30)

%%%%%%%%%%%%%%%%%%%%%%%%%%%%%%%%%%%%%%%%%%%%%%%%%%% FIRST RIGHT ELEMENT %%%%%%%%%%%%%%%%%%%%%%%%%%%%

\put(85.00,00.00){\makebox(0,0)[cc]{$\bullet$}}
\put(85.00,10.00){\makebox(0,0)[cc]{$\bullet$}}
\put(85.00,20.00){\makebox(0,0)[cc]{$\bullet$}}
\put(85.00,30.00){\makebox(0,0)[cc]{$\bullet$}}
\put(95.00,00.00){\makebox(0,0)[cc]{$\bullet$}}
\put(95.00,10.00){\makebox(0,0)[cc]{$\bullet$}}
\put(95.00,20.00){\makebox(0,0)[cc]{$\bullet$}}
\put(95.00,30.00){\makebox(0,0)[cc]{$\bullet$}}

\drawline(81.00,-04.00)(81.00,34.00)
\drawline(81.00,34.00)(99.00,34.00)
\drawline(99.00,34.00)(99.00,-04.00)
\drawline(99.00,-04.00)(81.00,-04.00)

\drawline(85.30,00.30)(95.30,00.30)
\drawline(85.30,10.30)(95.30,10.30)

\drawline(85.30,20.30)(87.30,20.30)
\drawline(87.30,20.30)(87.30,30.30)
\drawline(87.30,30.30)(85.30,30.30)

\drawline(95.30,20.30)(92.80,20.30)
\drawline(92.80,20.30)(92.80,30.30)
\drawline(92.80,30.30)(95.30,30.30)

%%%%%%%%%%%%%%%%%%%%%%%%%%%%%%%%%%%%%%%%%%%%%%%%%%% SECOND RIGHT ELEMENT %%%%%%%%%%%%%%%%%%%%%%%%%%%%

\put(108.00,00.00){\makebox(0,0)[cc]{$\bullet$}}
\put(108.00,10.00){\makebox(0,0)[cc]{$\bullet$}}
\put(108.00,20.00){\makebox(0,0)[cc]{$\bullet$}}
\put(108.00,30.00){\makebox(0,0)[cc]{$\bullet$}}
\put(118.00,00.00){\makebox(0,0)[cc]{$\bullet$}}
\put(118.00,10.00){\makebox(0,0)[cc]{$\bullet$}}
\put(118.00,20.00){\makebox(0,0)[cc]{$\bullet$}}
\put(118.00,30.00){\makebox(0,0)[cc]{$\bullet$}}

\drawline(104.00,-04.00)(104.00,34.00)
\drawline(104.00,34.00)(122.00,34.00)
\drawline(122.00,34.00)(122.00,-04.00)
\drawline(122.00,-04.00)(104.00,-04.00)

\drawline(108.30,10.30)(118.30,10.30)
\drawline(108.30,30.30)(118.30,30.30)

\drawline(108.30,00.30)(110.30,00.30)
\drawline(110.30,00.30)(110.30,20.30)
\drawline(110.30,20.30)(108.30,20.30)

\drawline(118.30,00.30)(115.80,00.30)
\drawline(115.80,00.30)(115.80,20.30)
\drawline(115.80,20.30)(118.30,20.30)

%%%%%%%%%%%%%%%%%%%%%%%%%%%%%%%%%%%%%%%%%%%%%%%%%%% THIRD RIGHT ELEMENT %%%%%%%%%%%%%%%%%%%%%%%%%%%%

\put(131.00,00.00){\makebox(0,0)[cc]{$\bullet$}}
\put(131.00,10.00){\makebox(0,0)[cc]{$\bullet$}}
\put(131.00,20.00){\makebox(0,0)[cc]{$\bullet$}}
\put(131.00,30.00){\makebox(0,0)[cc]{$\bullet$}}
\put(141.00,00.00){\makebox(0,0)[cc]{$\bullet$}}
\put(141.00,10.00){\makebox(0,0)[cc]{$\bullet$}}
\put(141.00,20.00){\makebox(0,0)[cc]{$\bullet$}}
\put(141.00,30.00){\makebox(0,0)[cc]{$\bullet$}}

\drawline(127.00,-04.00)(127.00,34.00)
\drawline(127.00,34.00)(145.00,34.00)
\drawline(145.00,34.00)(145.00,-04.00)
\drawline(145.00,-04.00)(127.00,-04.00)

\drawline(131.30,00.30)(141.30,00.30)
\drawline(131.30,20.30)(141.30,20.30)

\drawline(131.30,30.30)(133.30,30.30)
\drawline(133.30,30.30)(133.30,10.30)
\drawline(133.30,10.30)(131.30,10.30)

\drawline(141.30,30.30)(138.80,30.30)
\drawline(138.80,30.30)(138.80,10.30)
\drawline(138.80,10.30)(141.30,10.30)

%%%%%%%%%%%%%%%%%%%%%%%%%%%%%%%%%%%%%%%%%%%%%%%%%% 1-2 %%%%%%%%%%%%%%%%%%%%%%%%%%%%%%%%%%%%%%%%%%%%%%%

\drawline(21.30,00.30)(24.30,00.30)
\drawline(26.30,00.30)(29.30,00.30)
\drawline(31.30,00.30)(34.30,00.30)

\drawline(21.30,10.30)(24.30,10.30)
\drawline(26.30,10.30)(29.30,10.30)
\drawline(31.30,10.30)(34.30,10.30)

\drawline(21.30,20.30)(24.30,20.30)
\drawline(26.30,20.30)(29.30,20.30)
\drawline(31.30,20.30)(34.30,20.30)

\drawline(21.30,30.30)(24.30,30.30)
\drawline(26.30,30.30)(29.30,30.30)
\drawline(31.30,30.30)(34.30,30.30)

%%%%%%%%%%%%%%%%%%%%%%%%%%%%%%%%%%%%%%%%%%%%%%%%%% 2-3 %%%%%%%%%%%%%%%%%%%%%%%%%%%%%%%%%%%%%%%%%%%%%%%

\drawline(44.30,00.30)(47.30,00.30)
\drawline(49.30,00.30)(52.30,00.30)
\drawline(54.30,00.30)(57.30,00.30)

\drawline(44.30,10.30)(47.30,10.30)
\drawline(49.30,10.30)(52.30,10.30)
\drawline(54.30,10.30)(57.30,10.30)

\drawline(44.30,20.30)(47.30,20.30)
\drawline(49.30,20.30)(52.30,20.30)
\drawline(54.30,20.30)(57.30,20.30)

\drawline(44.30,30.30)(47.30,30.30)
\drawline(49.30,30.30)(52.30,30.30)
\drawline(54.30,30.30)(57.30,30.30)

%%%%%%%%%%%%%%%%%%%%%%%%%%%%%%%%%%%%%%%%%%%%%%%%%% 4-5 %%%%%%%%%%%%%%%%%%%%%%%%%%%%%%%%%%%%%%%%%%%%%%%

\drawline(95.30,00.30)(98.30,00.30)
\drawline(100.30,00.30)(103.30,00.30)
\drawline(105.30,00.30)(108.30,00.30)

\drawline(95.30,10.30)(98.30,10.30)
\drawline(100.30,10.30)(103.30,10.30)
\drawline(105.30,10.30)(108.30,10.30)

\drawline(95.30,20.30)(98.30,20.30)
\drawline(100.30,20.30)(103.30,20.30)
\drawline(105.30,20.30)(108.30,20.30)

\drawline(95.30,30.30)(98.30,30.30)
\drawline(100.30,30.30)(103.30,30.30)
\drawline(105.30,30.30)(108.30,30.30)

%%%%%%%%%%%%%%%%%%%%%%%%%%%%%%%%%%%%%%%%%%%%%%%%%% 5-6 %%%%%%%%%%%%%%%%%%%%%%%%%%%%%%%%%%%%%%%%%%%%%%%

\drawline(118.30,00.30)(121.30,00.30)
\drawline(123.30,00.30)(126.30,00.30)
\drawline(128.30,00.30)(131.30,00.30)

\drawline(118.30,10.30)(121.30,10.30)
\drawline(123.30,10.30)(126.30,10.30)
\drawline(128.30,10.30)(131.30,10.30)

\drawline(118.30,20.30)(121.30,20.30)
\drawline(123.30,20.30)(126.30,20.30)
\drawline(128.30,20.30)(131.30,20.30)

\drawline(118.30,30.30)(121.30,30.30)
\drawline(123.30,30.30)(126.30,30.30)
\drawline(128.30,30.30)(131.30,30.30)

%%%%%%%%%%%%%%%%%%%%%%%%%%%%%%%%%%%%%%%%%%%%%%%%%%%%%%%%%%%%%%%%%%%%%%%%%%%%%%%%%%%%%%%%%%%%%%%%%%

%%%%%%%%%%%%%%%%%%%%%%%%%%%%%%%%%%%%%%%%%%%%%%%%%%%%%%%%%%%%%%%%%%%%%%%%%%%%%%%%%%%%%%%%%
\put(76.00,15.00){\makebox(0,0)[cc]{$=$}}
\put(04.30,00.30){\makebox(0,0)[cc]{$l$}}
\put(04.30,10.30){\makebox(0,0)[cc]{$k$}}
\put(04.30,20.30){\makebox(0,0)[cc]{$j$}}
\put(04.30,30.30){\makebox(0,0)[cc]{$i$}}
\end{picture}
\caption{An example illustrating the  
relation~\eqref{6}.}\label{fig:example-for-6}

%%%%%%%%%%%%%%%%%%%%%%%%%%%%%%%%%%%%%%%%%%%%%%%%%%%%%%%%%%%%%%%%%%%%%%%%%%%%%%%%
%%%%
\end{figure}
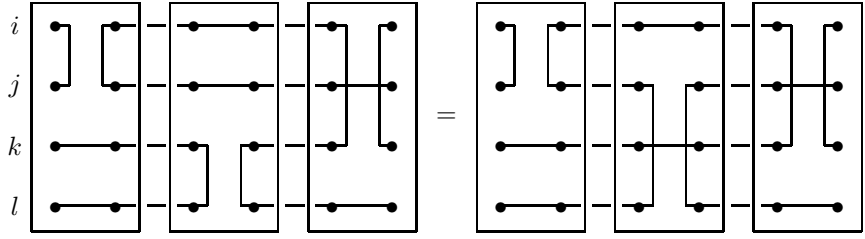

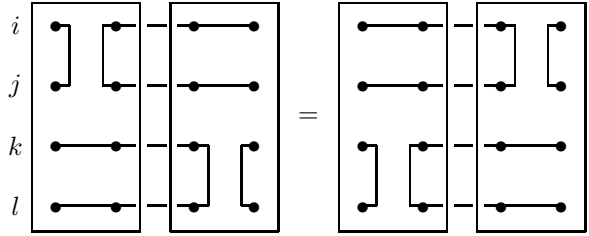
\begin{figure}
%%%%%%%%%%%%%%%%%%%%%%%%%%%%%%%%%%%%%%%%%%%%%%%%%%%%%%%%%%%%%%%%%%%%%%%%%%%%%%%%
%%%%%
%\input{figa1.pic}
%TexCad Options
%\grade{\on}
%\paths{\on}
%\beziermacro{\off}
%\reduce{\on}
%\snapping{\off}
%\quality{2.00}
%\graddiff{0.01}
%\snapasp{1}
%\zoom{1.00}
\special{em:linewidth 0.4pt} \unitlength 0.80mm
\linethickness{1pt}
\begin{picture}(150.00,35.00)
%%%%%%%%%%%%%%%%%%%%%%%%%%%%%%%%%%%%%%%%%%%%%%%%%%% SECOND LEFT ELEMENT %%%%%%%%%%%%%%%%%%%%%%%%%%%%

\put(34.00,00.00){\makebox(0,0)[cc]{$\bullet$}}
\put(34.00,10.00){\makebox(0,0)[cc]{$\bullet$}}
\put(34.00,20.00){\makebox(0,0)[cc]{$\bullet$}}
\put(34.00,30.00){\makebox(0,0)[cc]{$\bullet$}}
\put(44.00,00.00){\makebox(0,0)[cc]{$\bullet$}}
\put(44.00,10.00){\makebox(0,0)[cc]{$\bullet$}}
\put(44.00,20.00){\makebox(0,0)[cc]{$\bullet$}}
\put(44.00,30.00){\makebox(0,0)[cc]{$\bullet$}}

\drawline(30.00,-04.00)(30.00,34.00)
\drawline(30.00,34.00)(48.00,34.00)
\drawline(48.00,34.00)(48.00,-04.00)
\drawline(48.00,-04.00)(30.00,-04.00)

\drawline(34.30,00.30)(44.30,00.30)
\drawline(34.10,10.30)(44.30,10.30)

\drawline(34.30,20.30)(36.30,20.30)
\drawline(36.30,20.30)(36.30,30.30)
\drawline(36.30,30.30)(34.30,30.30)

\drawline(44.30,20.30)(41.80,20.30)
\drawline(41.80,20.30)(41.80,30.30)
\drawline(41.80,30.30)(44.30,30.30)

%%%%%%%%%%%%%%%%%%%%%%%%%%%%%%%%%%%%%%%%%%%%%%%%%%% THIRD LEFT ELEMENT %%%%%%%%%%%%%%%%%%%%%%%%%%%%

\put(57.00,00.00){\makebox(0,0)[cc]{$\bullet$}}
\put(57.00,10.00){\makebox(0,0)[cc]{$\bullet$}}
\put(57.00,20.00){\makebox(0,0)[cc]{$\bullet$}}
\put(57.00,30.00){\makebox(0,0)[cc]{$\bullet$}}
\put(67.00,00.00){\makebox(0,0)[cc]{$\bullet$}}
\put(67.00,10.00){\makebox(0,0)[cc]{$\bullet$}}
\put(67.00,20.00){\makebox(0,0)[cc]{$\bullet$}}
\put(67.00,30.00){\makebox(0,0)[cc]{$\bullet$}}

\drawline(53.00,-04.00)(53.00,34.00)
\drawline(53.00,34.00)(71.00,34.00)
\drawline(71.00,34.00)(71.00,-04.00)
\drawline(71.00,-04.00)(53.00,-04.00)

\drawline(57.30,20.30)(67.30,20.30)
\drawline(57.30,30.30)(67.30,30.30)

\drawline(57.30,00.30)(59.30,00.30)
\drawline(59.30,00.30)(59.30,10.30)
\drawline(59.30,10.30)(57.30,10.30)

\drawline(67.30,00.30)(64.80,00.30)
\drawline(64.80,00.30)(64.80,10.30)
\drawline(64.80,10.30)(67.30,10.30)

%%%%%%%%%%%%%%%%%%%%%%%%%%%%%%%%%%%%%%%%%%%%%%%%%%% FIRST RIGHT ELEMENT %%%%%%%%%%%%%%%%%%%%%%%%%%%%

\put(85.00,00.00){\makebox(0,0)[cc]{$\bullet$}}
\put(85.00,10.00){\makebox(0,0)[cc]{$\bullet$}}
\put(85.00,20.00){\makebox(0,0)[cc]{$\bullet$}}
\put(85.00,30.00){\makebox(0,0)[cc]{$\bullet$}}
\put(95.00,00.00){\makebox(0,0)[cc]{$\bullet$}}
\put(95.00,10.00){\makebox(0,0)[cc]{$\bullet$}}
\put(95.00,20.00){\makebox(0,0)[cc]{$\bullet$}}
\put(95.00,30.00){\makebox(0,0)[cc]{$\bullet$}}

\drawline(81.00,-04.00)(81.00,34.00)
\drawline(81.00,34.00)(99.00,34.00)
\drawline(99.00,34.00)(99.00,-04.00)
\drawline(99.00,-04.00)(81.00,-04.00)

\drawline(85.30,20.30)(95.30,20.30)
\drawline(85.30,30.30)(95.30,30.30)

\drawline(85.30,00.30)(87.30,00.30)
\drawline(87.30,00.30)(87.30,10.30)
\drawline(87.30,10.30)(85.30,10.30)

\drawline(95.30,00.30)(92.80,00.30)
\drawline(92.80,00.30)(92.80,10.30)
\drawline(92.80,10.30)(95.30,10.30)

%%%%%%%%%%%%%%%%%%%%%%%%%%%%%%%%%%%%%%%%%%%%%%%%%%% SECOND RIGHT ELEMENT %%%%%%%%%%%%%%%%%%%%%%%%%%%%

\put(108.00,00.00){\makebox(0,0)[cc]{$\bullet$}}
\put(108.00,10.00){\makebox(0,0)[cc]{$\bullet$}}
\put(108.00,20.00){\makebox(0,0)[cc]{$\bullet$}}
\put(108.00,30.00){\makebox(0,0)[cc]{$\bullet$}}
\put(118.00,00.00){\makebox(0,0)[cc]{$\bullet$}}
\put(118.00,10.00){\makebox(0,0)[cc]{$\bullet$}}
\put(118.00,20.00){\makebox(0,0)[cc]{$\bullet$}}
\put(118.00,30.00){\makebox(0,0)[cc]{$\bullet$}}

\drawline(104.00,-04.00)(104.00,34.00)
\drawline(104.00,34.00)(122.00,34.00)
\drawline(122.00,34.00)(122.00,-04.00)
\drawline(122.00,-04.00)(104.00,-04.00)

\drawline(108.30,10.30)(118.30,10.30)
\drawline(108.30,00.30)(118.30,00.30)

\drawline(108.30,30.30)(110.30,30.30)
\drawline(110.30,30.30)(110.30,20.30)
\drawline(110.30,20.30)(108.30,20.30)

\drawline(118.30,30.30)(115.80,30.30)
\drawline(115.80,30.30)(115.80,20.30)
\drawline(115.80,20.30)(118.30,20.30)

%%%%%%%%%%%%%%%%%%%%%%%%%%%%%%%%%%%%%%%%%%%%%%%%%% 2-3 %%%%%%%%%%%%%%%%%%%%%%%%%%%%%%%%%%%%%%%%%%%%%%%

\drawline(44.30,00.30)(47.30,00.30)
\drawline(49.30,00.30)(52.30,00.30)
\drawline(54.30,00.30)(57.30,00.30)

\drawline(44.30,10.30)(47.30,10.30)
\drawline(49.30,10.30)(52.30,10.30)
\drawline(54.30,10.30)(57.30,10.30)

\drawline(44.30,20.30)(47.30,20.30)
\drawline(49.30,20.30)(52.30,20.30)
\drawline(54.30,20.30)(57.30,20.30)

\drawline(44.30,30.30)(47.30,30.30)
\drawline(49.30,30.30)(52.30,30.30)
\drawline(54.30,30.30)(57.30,30.30)

%%%%%%%%%%%%%%%%%%%%%%%%%%%%%%%%%%%%%%%%%%%%%%%%%% 4-5 %%%%%%%%%%%%%%%%%%%%%%%%%%%%%%%%%%%%%%%%%%%%%%%

\drawline(95.30,00.30)(98.30,00.30)
\drawline(100.30,00.30)(103.30,00.30)
\drawline(105.30,00.30)(108.30,00.30)

\drawline(95.30,10.30)(98.30,10.30)
\drawline(100.30,10.30)(103.30,10.30)
\drawline(105.30,10.30)(108.30,10.30)

\drawline(95.30,20.30)(98.30,20.30)
\drawline(100.30,20.30)(103.30,20.30)
\drawline(105.30,20.30)(108.30,20.30)

\drawline(95.30,30.30)(98.30,30.30)
\drawline(100.30,30.30)(103.30,30.30)
\drawline(105.30,30.30)(108.30,30.30)

%%%%%%%%%%%%%%%%%%%%%%%%%%%%%%%%%%%%%%%%%%%%%%%%%%%%%%%%%%%%%%%%%%%%%%%%%%%%%%%%%%%%%%%%%
\put(76.00,15.00){\makebox(0,0)[cc]{$=$}}
\put(27.30,00.30){\makebox(0,0)[cc]{$l$}}
\put(27.30,10.30){\makebox(0,0)[cc]{$k$}}
\put(27.30,20.30){\makebox(0,0)[cc]{$j$}}
\put(27.30,30.30){\makebox(0,0)[cc]{$i$}}
\end{picture}
\caption{An example illustrating the
 relation~\eqref{7}.}\label{fig:example-for-7}

%%%%%%%%%%%%%%%%%%%%%%%%%%%%%%%%%%%%%%%%%%%%%%%%%%%%%%%%%%%%%%%%%%%%%%%%%%%%%%%%
%%%%
\end{figure}

A straightforward calculation shows that the generators
$\sigma_{i,j}$ of $\B_n\setminus\Sym_n$ satisfy the relations
\eqref{1}--\eqref{7} (the relations \eqref{1} and \eqref{2} are obvious,
and the relations \eqref{3}--\eqref{7} are illustrated on Figure~\ref{fig:example-for-3}, Figure~\ref{fig:example-for-4},
Figure~\ref{fig:example-for-5}, Figure~\ref{fig:example-for-6} 
and Figure~\ref{fig:example-for-7}). Thus there is a homomorphism $\varphi:T\to\B_n\setminus\Sym_n$, sending $\tau_{i,j}$ to $\sigma_{i,j}$. 
Our main goal in the section is to prove the following theorem:

%%%%%%%%%%%%%%%%%%%%%%%%%%%%%%%%%%%%%%%%%%%%%%%%%%%%%%%%%%%%%%%%%%%%%%%%%%%%%%%%%%%%%%%%%%%%%%%%%%%%%%%

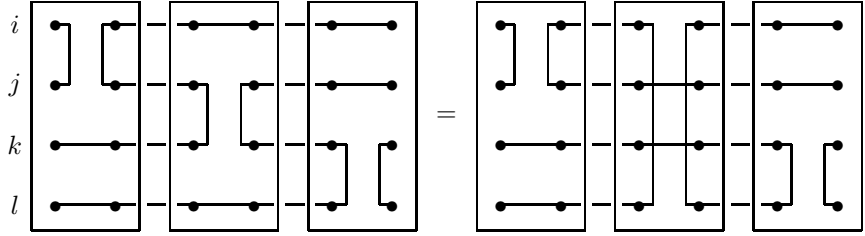
\begin{figure}
%%%%%%%%%%%%%%%%%%%%%%%%%%%%%%%%%%%%%%%%%%%%%%%%%%%%%%%%%%%%%%%%%%%%%%%%%%%%%%%%
%%%%%
%\input{figa1.pic}
%TexCad Options
%\grade{\on}
%\paths{\on}
%\beziermacro{\off}
%\reduce{\on}
%\snapping{\off}
%\quality{2.00}
%\graddiff{0.01}
%\snapasp{1}
%\zoom{1.00}
\special{em:linewidth 0.4pt} \unitlength 0.80mm
\linethickness{1pt}
\begin{picture}(150.00,35.00)
%%%%%%%%%%%%%%%%%%%%%%%%%%%%%%%%%%%%%%%%%%%%%%%% FIRST LEFT ELEMENT %%%%%%%%%%%%%%%%%%%%%%%%%%%%%%%%%%%
\put(11.00,00.00){\makebox(0,0)[cc]{$\bullet$}}
\put(11.00,10.00){\makebox(0,0)[cc]{$\bullet$}}
\put(11.00,20.00){\makebox(0,0)[cc]{$\bullet$}}
\put(11.00,30.00){\makebox(0,0)[cc]{$\bullet$}}
\put(21.00,00.00){\makebox(0,0)[cc]{$\bullet$}}
\put(21.00,10.00){\makebox(0,0)[cc]{$\bullet$}}
\put(21.00,20.00){\makebox(0,0)[cc]{$\bullet$}}
\put(21.00,30.00){\makebox(0,0)[cc]{$\bullet$}}

\drawline(07.00,-04.00)(07.00,34.00)
\drawline(07.00,34.00)(25.00,34.00)
\drawline(25.00,34.00)(25.00,-04.00)
\drawline(25.00,-04.00)(07.00,-04.00)

\drawline(11.30,00.30)(21.30,00.30)
\drawline(11.30,10.30)(21.30,10.30)

\drawline(11.30,20.30)(13.30,20.30)
\drawline(13.30,20.30)(13.30,30.30)
\drawline(13.30,30.30)(11.30,30.30)

\drawline(21.30,20.30)(18.80,20.30)
\drawline(18.80,20.30)(18.80,30.30)
\drawline(18.80,30.30)(21.30,30.30)

%%%%%%%%%%%%%%%%%%%%%%%%%%%%%%%%%%%%%%%%%%%%%%%%%%% SECOND LEFT ELEMENT %%%%%%%%%%%%%%%%%%%%%%%%%%%%

\put(34.00,00.00){\makebox(0,0)[cc]{$\bullet$}}
\put(34.00,10.00){\makebox(0,0)[cc]{$\bullet$}}
\put(34.00,20.00){\makebox(0,0)[cc]{$\bullet$}}
\put(34.00,30.00){\makebox(0,0)[cc]{$\bullet$}}
\put(44.00,00.00){\makebox(0,0)[cc]{$\bullet$}}
\put(44.00,10.00){\makebox(0,0)[cc]{$\bullet$}}
\put(44.00,20.00){\makebox(0,0)[cc]{$\bullet$}}
\put(44.00,30.00){\makebox(0,0)[cc]{$\bullet$}}

\drawline(30.00,-04.00)(30.00,34.00)
\drawline(30.00,34.00)(48.00,34.00)
\drawline(48.00,34.00)(48.00,-04.00)
\drawline(48.00,-04.00)(30.00,-04.00)

\drawline(34.30,00.30)(44.30,00.30)
\drawline(34.30,30.30)(44.30,30.30)

\drawline(34.30,10.30)(36.30,10.30)
\drawline(36.30,10.30)(36.30,20.30)
\drawline(36.30,20.30)(34.30,20.30)

\drawline(44.30,10.30)(41.80,10.30)
\drawline(41.80,10.30)(41.80,20.30)
\drawline(41.80,20.30)(44.30,20.30)

%%%%%%%%%%%%%%%%%%%%%%%%%%%%%%%%%%%%%%%%%%%%%%%%%%% THIRD LEFT ELEMENT %%%%%%%%%%%%%%%%%%%%%%%%%%%%

\put(57.00,00.00){\makebox(0,0)[cc]{$\bullet$}}
\put(57.00,10.00){\makebox(0,0)[cc]{$\bullet$}}
\put(57.00,20.00){\makebox(0,0)[cc]{$\bullet$}}
\put(57.00,30.00){\makebox(0,0)[cc]{$\bullet$}}
\put(67.00,00.00){\makebox(0,0)[cc]{$\bullet$}}
\put(67.00,10.00){\makebox(0,0)[cc]{$\bullet$}}
\put(67.00,20.00){\makebox(0,0)[cc]{$\bullet$}}
\put(67.00,30.00){\makebox(0,0)[cc]{$\bullet$}}

\drawline(53.00,-04.00)(53.00,34.00)
\drawline(53.00,34.00)(71.00,34.00)
\drawline(71.00,34.00)(71.00,-04.00)
\drawline(71.00,-04.00)(53.00,-04.00)

\drawline(57.30,20.30)(67.30,20.30)
\drawline(57.30,30.30)(67.30,30.30)

\drawline(57.30,00.30)(59.30,00.30)
\drawline(59.30,00.30)(59.30,10.30)
\drawline(59.30,10.30)(57.30,10.30)

\drawline(67.30,00.30)(64.80,00.30)
\drawline(64.80,00.30)(64.80,10.30)
\drawline(64.80,10.30)(67.30,10.30)

%%%%%%%%%%%%%%%%%%%%%%%%%%%%%%%%%%%%%%%%%%%%%%%%%%% FIRST RIGHT ELEMENT %%%%%%%%%%%%%%%%%%%%%%%%%%%%

\put(85.00,00.00){\makebox(0,0)[cc]{$\bullet$}}
\put(85.00,10.00){\makebox(0,0)[cc]{$\bullet$}}
\put(85.00,20.00){\makebox(0,0)[cc]{$\bullet$}}
\put(85.00,30.00){\makebox(0,0)[cc]{$\bullet$}}
\put(95.00,00.00){\makebox(0,0)[cc]{$\bullet$}}
\put(95.00,10.00){\makebox(0,0)[cc]{$\bullet$}}
\put(95.00,20.00){\makebox(0,0)[cc]{$\bullet$}}
\put(95.00,30.00){\makebox(0,0)[cc]{$\bullet$}}

\drawline(81.00,-04.00)(81.00,34.00)
\drawline(81.00,34.00)(99.00,34.00)
\drawline(99.00,34.00)(99.00,-04.00)
\drawline(99.00,-04.00)(81.00,-04.00)

\drawline(85.30,00.30)(95.30,00.30)
\drawline(85.30,10.30)(95.30,10.30)

\drawline(85.30,20.30)(87.30,20.30)
\drawline(87.30,20.30)(87.30,30.30)
\drawline(87.30,30.30)(85.30,30.30)

\drawline(95.30,20.30)(92.80,20.30)
\drawline(92.80,20.30)(92.80,30.30)
\drawline(92.80,30.30)(95.30,30.30)

%%%%%%%%%%%%%%%%%%%%%%%%%%%%%%%%%%%%%%%%%%%%%%%%%%% SECOND RIGHT ELEMENT %%%%%%%%%%%%%%%%%%%%%%%%%%%%

\put(108.00,00.00){\makebox(0,0)[cc]{$\bullet$}}
\put(108.00,10.00){\makebox(0,0)[cc]{$\bullet$}}
\put(108.00,20.00){\makebox(0,0)[cc]{$\bullet$}}
\put(108.00,30.00){\makebox(0,0)[cc]{$\bullet$}}
\put(118.00,00.00){\makebox(0,0)[cc]{$\bullet$}}
\put(118.00,10.00){\makebox(0,0)[cc]{$\bullet$}}
\put(118.00,20.00){\makebox(0,0)[cc]{$\bullet$}}
\put(118.00,30.00){\makebox(0,0)[cc]{$\bullet$}}

\drawline(104.00,-04.00)(104.00,34.00)
\drawline(104.00,34.00)(122.00,34.00)
\drawline(122.00,34.00)(122.00,-04.00)
\drawline(122.00,-04.00)(104.00,-04.00)

\drawline(108.30,10.30)(118.30,10.30)
\drawline(108.30,20.30)(118.30,20.30)

\drawline(108.30,00.30)(110.30,00.30)
\drawline(110.30,00.30)(110.30,30.30)
\drawline(110.30,30.30)(108.30,30.30)

\drawline(118.30,00.30)(115.80,00.30)
\drawline(115.80,00.30)(115.80,30.30)
\drawline(115.80,30.30)(118.30,30.30)

%%%%%%%%%%%%%%%%%%%%%%%%%%%%%%%%%%%%%%%%%%%%%%%%%%% THIRD RIGHT ELEMENT %%%%%%%%%%%%%%%%%%%%%%%%%%%%

\put(131.00,00.00){\makebox(0,0)[cc]{$\bullet$}}
\put(131.00,10.00){\makebox(0,0)[cc]{$\bullet$}}
\put(131.00,20.00){\makebox(0,0)[cc]{$\bullet$}}
\put(131.00,30.00){\makebox(0,0)[cc]{$\bullet$}}
\put(141.00,00.00){\makebox(0,0)[cc]{$\bullet$}}
\put(141.00,10.00){\makebox(0,0)[cc]{$\bullet$}}
\put(141.00,20.00){\makebox(0,0)[cc]{$\bullet$}}
\put(141.00,30.00){\makebox(0,0)[cc]{$\bullet$}}

\drawline(127.00,-04.00)(127.00,34.00)
\drawline(127.00,34.00)(145.00,34.00)
\drawline(145.00,34.00)(145.00,-04.00)
\drawline(145.00,-04.00)(127.00,-04.00)

\drawline(131.30,30.30)(141.30,30.30)
\drawline(131.30,20.30)(141.30,20.30)

\drawline(131.30,00.30)(133.30,00.30)
\drawline(133.30,00.30)(133.30,10.30)
\drawline(133.30,10.30)(131.30,10.30)

\drawline(141.30,00.30)(138.80,00.30)
\drawline(138.80,00.30)(138.80,10.30)
\drawline(138.80,10.30)(141.30,10.30)

%%%%%%%%%%%%%%%%%%%%%%%%%%%%%%%%%%%%%%%%%%%%%%%%%% 1-2 %%%%%%%%%%%%%%%%%%%%%%%%%%%%%%%%%%%%%%%%%%%%%%%

\drawline(21.30,00.30)(24.30,00.30)
\drawline(26.30,00.30)(29.30,00.30)
\drawline(31.30,00.30)(34.30,00.30)

\drawline(21.30,10.30)(24.30,10.30)
\drawline(26.30,10.30)(29.30,10.30)
\drawline(31.30,10.30)(34.30,10.30)

\drawline(21.30,20.30)(24.30,20.30)
\drawline(26.30,20.30)(29.30,20.30)
\drawline(31.30,20.30)(34.30,20.30)

\drawline(21.30,30.30)(24.30,30.30)
\drawline(26.30,30.30)(29.30,30.30)
\drawline(31.30,30.30)(34.30,30.30)

%%%%%%%%%%%%%%%%%%%%%%%%%%%%%%%%%%%%%%%%%%%%%%%%%% 2-3 %%%%%%%%%%%%%%%%%%%%%%%%%%%%%%%%%%%%%%%%%%%%%%%

\drawline(44.30,00.30)(47.30,00.30)
\drawline(49.30,00.30)(52.30,00.30)
\drawline(54.30,00.30)(57.30,00.30)

\drawline(44.30,10.30)(47.30,10.30)
\drawline(49.30,10.30)(52.30,10.30)
\drawline(54.30,10.30)(57.30,10.30)

\drawline(44.30,20.30)(47.30,20.30)
\drawline(49.30,20.30)(52.30,20.30)
\drawline(54.30,20.30)(57.30,20.30)

\drawline(44.30,30.30)(47.30,30.30)
\drawline(49.30,30.30)(52.30,30.30)
\drawline(54.30,30.30)(57.30,30.30)

%%%%%%%%%%%%%%%%%%%%%%%%%%%%%%%%%%%%%%%%%%%%%%%%%% 4-5 %%%%%%%%%%%%%%%%%%%%%%%%%%%%%%%%%%%%%%%%%%%%%%%

\drawline(95.30,00.30)(98.30,00.30)
\drawline(100.30,00.30)(103.30,00.30)
\drawline(105.30,00.30)(108.30,00.30)

\drawline(95.30,10.30)(98.30,10.30)
\drawline(100.30,10.30)(103.30,10.30)
\drawline(105.30,10.30)(108.30,10.30)

\drawline(95.30,20.30)(98.30,20.30)
\drawline(100.30,20.30)(103.30,20.30)
\drawline(105.30,20.30)(108.30,20.30)

\drawline(95.30,30.30)(98.30,30.30)
\drawline(100.30,30.30)(103.30,30.30)
\drawline(105.30,30.30)(108.30,30.30)

%%%%%%%%%%%%%%%%%%%%%%%%%%%%%%%%%%%%%%%%%%%%%%%%%% 5-6 %%%%%%%%%%%%%%%%%%%%%%%%%%%%%%%%%%%%%%%%%%%%%%%

\drawline(118.30,00.30)(121.30,00.30)
\drawline(123.30,00.30)(126.30,00.30)
\drawline(128.30,00.30)(131.30,00.30)

\drawline(118.30,10.30)(121.30,10.30)
\drawline(123.30,10.30)(126.30,10.30)
\drawline(128.30,10.30)(131.30,10.30)

\drawline(118.30,20.30)(121.30,20.30)
\drawline(123.30,20.30)(126.30,20.30)
\drawline(128.30,20.30)(131.30,20.30)

\drawline(118.30,30.30)(121.30,30.30)
\drawline(123.30,30.30)(126.30,30.30)
\drawline(128.30,30.30)(131.30,30.30)

%%%%%%%%%%%%%%%%%%%%%%%%%%%%%%%%%%%%%%%%%%%%%%%%%%%%%%%%%%%%%%%%%%%%%%%%%%%%%%%%%%%%%%%%%%%%%%%%%%

%%%%%%%%%%%%%%%%%%%%%%%%%%%%%%%%%%%%%%%%%%%%%%%%%%%%%%%%%%%%%%%%%%%%%%%%%%%%%%%%%%%%%%%%%
\put(76.00,15.00){\makebox(0,0)[cc]{$=$}}

\put(04.30,00.30){\makebox(0,0)[cc]{$l$}}
\put(04.30,10.30){\makebox(0,0)[cc]{$k$}}
\put(04.30,20.30){\makebox(0,0)[cc]{$j$}}
\put(04.30,30.30){\makebox(0,0)[cc]{$i$}}
\end{picture}
\caption{An example illustrating the 
 relation~\eqref{3}.}\label{fig:example-for-3}

%%%%%%%%%%%%%%%%%%%%%%%%%%%%%%%%%%%%%%%%%%%%%%%%%%%%%%%%%%%%%%%%%%%%%%%%%%%%%%%%
%%%%
\end{figure}

%%%%%%%%%%%%%%%%%%%%%%%%%%%%%%%%%%%%%%%%%%%%%%%%%%%%%%%%%%%%%%%%%%%%%%%%%%%%%%%%%%%%%%%%%%%%%%%%%%%%%%%

%%%%%%%%%%%%%%%%%%%%%%%%%%%%%%%%%%%%%%%%%%%%%%%%%%%%%%%%%%%%%%%%%%%%%%%%%%%%%%%%%%%%%%%%%%%%%%%%%%%%%%%

%%%%%%%%%%%%%%%%%%%%%%%%%%%%%%%%%%%%%%%%%%%%%%%%%%%%%%%%%%%%%%%%%%%%%%%%%%%%%%%%%%%%%%%%%%%%%%%%%%%%%%%

%%%%%%%%%%%%%%%%%%%%%%%%%%%%%%%%%%%%%%%%%%%%%%%%%%%%%%%%%%%%%%%%%%%%%%%%%%%%%%%%%%%%%%%%%%%%%%%%%%%%%%%

\begin{theorem}\label{th:main}
$\varphi:T\to \B_n\setminus \Sym_n$ is an isomorphism. 
\end{theorem}

The rest of this section is devoted to the proof of Theorem~\ref{th:main},
which we will divide into steps, formulated as lemmas and propositions.
To distinguish $\tau_{i,j}$ from the atoms $\sigma_{i,j}$ we will call
$\tau_{i,j}$ \emph{quarks}. Two quarks $\tau_{i,j}$ and $\tau_{k,l}$ 
are said to be {\em connected} provided that 
$\{i,j\}\cap\{k,l\}\ne\varnothing$. We denote by $\A=\A_n$ the set of 
all quarks (the alphabet of our presentation for $T$), and by $\A^{+}$ 
the  free semigroup over $\A$. In what follows we will do all our 
computations with words in $T$, not $\A^{+}$. In particular, $v=w$ for
$v,w\in \A^+$ means that $v=w$ in $T$. 

A word, $\tau_{i_1,j_1}\tau_{i_2,j_2}\dots\tau_{i_k,j_k}\in\A^{+}$, 
will be called {\em connected} if $\tau_{i_s,j_s}$ and 
$\tau_{i_{s+1},j_{s+1}}$ are connected for all $1\leq s\leq k-1$.
We start with the following statement:

\begin{proposition}\label{p1}
Each element of the semigroup $T$ can we written in the form
$w\tau_{i_1,j_1}\tau_{i_2,j_2}\dots \tau_{i_k,j_k}$, where 
$w\tau_{i_1,j_1}\in \A^+$ is connected and all sets 
$\{i_s,j_s\}$, $s=1,\dots,k$, are pairwise disjoint.
\end{proposition}

\begin{proof}
We use induction on the length of the element. For elements of 
length $1$ the statement is obvious. Let 
$v=w\tau_{i_1,j_1}\tau_{i_2,j_2}\dots \tau_{i_k,j_k}\in T$ be such 
that $w\tau_{i_1,j_1}\in \A^+$ is connected and all sets 
$\{i_s,j_s\}$, $s=1,\dots,k$, are pairwise disjoint. Let further
$\tau_{i,j}$ be some generator. To complete the proof we have to 
show that the element $v\tau_{i,j}$ can be written in the necessary 
form. Without loss of generality we can assume that we have one
of the following cases:

{\bf Case~1:} the set $\{i,j\}$ is disjoint with all 
$\{i_s,j_s\}$, $s=1,\dots,k$. In this case the statement is trivial. 

{\bf Case~2:} the set $\{i,j\}$ is disjoint with all 
$\{i_s,j_s\}$, $s=2,\dots,k$, but not with $\tau_{i_1,j_1}$. In this case
we can use \eqref{7} to write 
\begin{displaymath}
v\tau_{i,j}=w\tau_{i_1,j_1}\tau_{i,j}\tau_{i_2,j_2}\dots \tau_{i_k,j_k}.
\end{displaymath}
Observe that $w\tau_{i_1,j_1}\tau_{i,j}$ is connected, and the necessary
statement follows again.

{\bf Case~3:} $i=i_1$ and $j\in \cup_{s=2}^k\{i_s,j_s\}$. Using
\eqref{7} we can even assume $j=j_2$. Using \eqref{7} and \eqref{6} we have
\begin{displaymath}
v\tau_{i,j}=w\tau_{i,j_1}\tau_{i_2,j}\tau_{i,j}
\tau_{i_3,j_3}\dots \tau_{i_k,j_k}=
w\tau_{i,j_1}\tau_{i_2,j_1}\tau_{i,j}
\tau_{i_3,j_3}\dots \tau_{i_k,j_k}.
\end{displaymath}
Here $w\tau_{i,j_1}\tau_{i_2,j_1}$ is connected and 
the sets $\{i_2,j_1\}$, $\{i,j\}$, $\{i_s,j_s\}$, $s=3,\dots,k$, are
disjoint. The claim follows.

{\bf Case~4:} $i\in \cup_{s=2}^k\{i_s,j_s\}$ and 
$j\not\in \cup_{s=1}^k\{i_s,j_s\}$. Using
\eqref{7},  we can even assume $i=i_2$. 
In this case we can use \eqref{7} to write 
\begin{equation}\label{a1}
v\tau_{i,j}=w\tau_{i_1,j_1}\tau_{i,j_2}\tau_{i,j}
\tau_{i_3,j_3}\dots \tau_{i_k,j_k}.
\end{equation}
Now we have:
\begin{equation}\label{a2}
\begin{array}{rclc}
\tau_{i_1,j_1}\tau_{i,j_2}\tau_{i,j} & = & 
\tau_{i_1,j_1}\tau_{i,j_2}\tau_{i,j}\tau_{i_1,j}\tau_{i,j}&
\text{(by \eqref{5})}\\
& = & 
\tau_{i_1,j_1}\tau_{i,j_2}\tau_{i_1,j_2}\tau_{i_1,j}\tau_{i,j}&
\text{(by \eqref{3})}\\
& = & 
\tau_{i_1,j_1}\tau_{i,j_1}\tau_{i_1,j_2}\tau_{i_1,j}\tau_{i,j}&
\text{(by \eqref{6})}\\
& = & 
\tau_{i_1,j_1}\tau_{i_1,j_2}\tau_{i_1,j}\tau_{i,j_1}\tau_{i,j}&
\text{(by \eqref{7})}\\
& = & 
\tau_{i_1,j_1}\tau_{i_1,j_2}\tau_{i_1,j}\tau_{i_1,j_1}\tau_{i,j}&
\text{(by \eqref{6})}.
\end{array}
\end{equation}
From \eqref{a1} and \eqref{a2} we have:
\begin{displaymath}
v\tau_{i,j}= w\tau_{i_1,j_1}\tau_{i,j_2}\tau_{i,j}
\tau_{i_3,j_3}\dots \tau_{i_k,j_k}=
w\tau_{i_1,j_1}\tau_{i_1,j_2}\tau_{i_1,j}\tau_{i_1,j_1}\tau_{i,j}
\tau_{i_3,j_3}\dots \tau_{i_k,j_k}.
\end{displaymath}
Here $w\tau_{i_1,j_1}\tau_{i_1,j_2}\tau_{i_1,j}\tau_{i_1,j_1}$ 
is connected and the sets $\{i_1,j_1\}$, $\{i,j\}$, 
$\{i_s,j_s\}$, $s=3,\dots,k$, are disjoint. The claim follows.

{\bf Case~5:} $i,j\in \cup_{s=2}^k\{i_s,j_s\}$. If 
$\{i,j\}=\{i_s,j_s\}$ for some $s\geq 2$, the statement follows
from \eqref{7} and \eqref{2}. Otherwise, using \eqref{7} we
can even assume $i=i_2$, $j=j_3$.
In this case we can use \eqref{7} to write 
\begin{equation}\label{a3}
v\tau_{i,j}=w\tau_{i_1,j_1}\tau_{i,j_2}\tau_{i_3,j}
\tau_{i,j}\tau_{i_4,j_4}\dots \tau_{i_k,j_k}.
\end{equation}
Now we have:
\begin{equation}\label{a4}
\begin{array}{rclc}
\tau_{i_1,j_1}\tau_{i,j_2}\tau_{i_3,j}\tau_{i,j} & = & 
\tau_{i_1,j_1}\tau_{i,j_2}\tau_{i_3,j_2}\tau_{i,j}&
\text{(by \eqref{6})}\\
& = & 
\tau_{i,j_2}\tau_{i_3,j_2}\tau_{i_1,j_1}\tau_{i,j}&
\text{(by \eqref{7})}\\
& = & 
\tau_{i,j_2}\tau_{i_3,j_2}\tau_{i_1,j_1}\tau_{j_2,j_1}\tau_{i_1,j_1}\tau_{i,j}&
\text{(by \eqref{5})}\\
& = & 
\tau_{i,j_2}\tau_{i_1,j_1}\tau_{j_2,i_3}\tau_{j_2,j_1}\tau_{i_1,j_1}\tau_{i,j}&
\text{(by \eqref{7})}\\
& = & 
\tau_{i,j_2}\tau_{i_1,j_1}\tau_{i_1,i_3}\tau_{j_2,j_1}\tau_{i_1,j_1}\tau_{i,j}&
\text{(by \eqref{6})}\\
& = & 
\tau_{i_1,j_1}\tau_{i_1,i_3}\tau_{i,j_2}\tau_{j_2,j_1}\tau_{i_1,j_1}\tau_{i,j}&
\text{(by \eqref{7})}\\
& = & 
\tau_{i_1,j_1}\tau_{i_1,i_3}\tau_{i,j_2}\tau_{i,i_1}\tau_{i_1,j_1}\tau_{i,j}&
\text{(by \eqref{3})}\\
& = & 
\tau_{i_1,j_1}\tau_{i_1,i_3}\tau_{i_3,j_2}\tau_{i,i_1}\tau_{i_1,j_1}\tau_{i,j}&
\text{(by \eqref{6})}\\
& = & 
\tau_{i_1,j_1}\tau_{i_1,i_3}\tau_{i,i_1}\tau_{i_1,j_1}\tau_{i_3,j_2}\tau_{i,j}&
\text{(by \eqref{7})}.
\end{array}
\end{equation}
From \eqref{a3} and \eqref{a4} we have:
\begin{multline*}
v\tau_{i,j}= w\tau_{i_1,j_1}\tau_{i,j_2}\tau_{i_3,j}
\tau_{i,j}\tau_{i_4,j_4}\dots \tau_{i_k,j_k}=\\
w\tau_{i_1,j_1}\tau_{i_1,i_3}\tau_{i,i_1}\tau_{i_1,j_1}
\tau_{i_3,j_2}\tau_{i,j} \tau_{i_4,j_4}\dots \tau_{i_k,j_k}.
\end{multline*}
Here $w\tau_{i_1,j_1}\tau_{i_1,i_3}\tau_{i,i_1}\tau_{i_1,j_1}$ 
is connected and the sets $\{i_1,j_1\}$, $\{i_3,j_2\}$, $\{i,j\}$, 
$\{i_s,j_s\}$, $s=4,\dots,k$, are disjoint. The claim follows.

Now the proof is completed by induction.
\end{proof}

\begin{lemma}\label{l2}
There is a unique anti-involution, $*:T\to T$, satisfying
$\tau_{i,j}^*=\tau_{i,j}$ for all $i,j\in\{1,2,\dots,n\}$, $i\neq j$.
\end{lemma}

\begin{proof}
Existence follows from the fact that the relations
\eqref{1}--\eqref{7} are stable with respect to
$*$. The uniqueness follows from the fact that 
$T$ is generated by $\tau_{i,j}$, $i\neq j\in\{1,2,\dots,n\}$.
\end{proof}

\begin{lemma}\label{l3}
Let $\tau_{i,j}w\in\A^+$ be connected. Then
$(\tau_{i,j}w)(\tau_{i,j}w)^*=\tau_{i,j}$.
\end{lemma}

\begin{proof}
Let $w=\tau_{i_1,j_1}\dots\tau_{i_k,j_k}$. Since
$\tau_{i,j}w$ is connected, applying \eqref{5}, \eqref{2} 
and the definition of $*$, we compute:
\begin{displaymath}
\begin{array}{rcl}
(\tau_{i,j}w)(\tau_{i,j}w)^* &= &
\tau_{i,j} \tau_{i_1,j_1}\dots\tau_{i_k,j_k}
\tau_{i_k,j_k}\dots\tau_{i_1,j_1}\tau_{i,j}\\
&= &
\tau_{i,j} \tau_{i_1,j_1}\dots\tau_{i_{k-1},j_{k-1}}
\tau_{i_k,j_k}\tau_{i_{k-1},j_{k-1}}\dots\tau_{i_1,j_1}\tau_{i,j}\\
&= &
\tau_{i,j} \tau_{i_1,j_1}\dots\tau_{i_{k-2},j_{k-2}}
\tau_{i_{k-1},j_{k-1}}\tau_{i_{k-2},j_{k-2}}\dots\tau_{i_1,j_1}\tau_{i,j}\\
&\dots &\\
&= &
\tau_{i,j} \tau_{i_1,j_1}\tau_{i,j}\\
&= &
\tau_{i,j}.
\end{array}
\end{displaymath}
\end{proof}

The elements of the form $\tau_{i_1,j_1}\dots\tau_{i_k,j_k}$,
where $\{i_s,j_s\}$, $s=1,\dots,k$, are pairwise disjoint, will
be called {\em standard idempotents}. That these
elements are indeed idempotents, follows immediately from
\eqref{7} and \eqref{2}.

\begin{corollary}\label{cgreen}
\begin{enumerate}[(i)]
\item \label{cgreen.1} Every element of $T$ is $\mathcal{L}$-equivalent
to a standard idempotent.
\item \label{cgreen.2} $T$ is regular.
\item \label{cgreen.3} The map $\varphi$ induces a bijection 
between the sets of $\mathcal{L}$-classes for the semigroups
$T$ and $\B_n\setminus\Sym_n$. Similar for the $\mathcal{R}$-, $\mathcal{H}$-, 
and $\mathcal{D}$-classes. 
\end{enumerate}
\end{corollary}

\begin{proof}
Let $v\in\A^+$. By Proposition~\ref{p1} we can write
$v=w\tau_{i_1,j_1}\tau_{i_2,j_2}\dots \tau_{i_k,j_k}$, where 
$w\tau_{i_1,j_1}\in \A^+$ is connected and all sets 
$\{i_s,j_s\}$, $s=1,\dots,k$, are pairwise disjoint.
By definition, the element 
$\epsilon=\tau_{i_1,j_1}\tau_{i_2,j_2}\dots \tau_{i_k,j_k}$
is standard. We obviously have $v=w\tau_{i_1,j_1}\epsilon$.
By Lemma~\ref{l3} we have
\begin{displaymath}
\tau_{i_1,j_1}w^*v=
\tau_{i_1,j_1}w^*w\tau_{i_1,j_1}\tau_{i_2,j_2}\dots \tau_{i_k,j_k}=
\tau_{i_1,j_1}\tau_{i_2,j_2}\dots \tau_{i_k,j_k}=\epsilon.
\end{displaymath}
Hence $v\mathcal{L}\epsilon$, which proves \eqref{cgreen.1}.
\eqref{cgreen.1} implies that every $\mathcal{L}$-class of
$T$ contains an idempotent, and hence \eqref{cgreen.2} follows.

By Lemma~\ref{lm:Green-on-Brauer}, the images of
standard idempotents under $\varphi$ belong to different
$\mathcal{L}$-classes of $\B_n\setminus\Sym_n$. Hence different
standard idempotents of $T$ belong to different 
$\mathcal{L}$-classes of $T$. In particular, there is a
bijection between $\mathcal{L}$-classes of $T$ and
standard idempotents. Since $\varphi$ is surjective, there is
also a bijection between $\mathcal{L}$-classes of $\B_n\setminus\Sym_n$ and
standard idempotents. This implies \eqref{cgreen.3} for 
$\mathcal{L}$-classes. For $\mathcal{R}$-classes the statement now
follows by applying $*$. For $\mathcal{H}$- and $\mathcal{D}$-classes 
the statement follows from the definition and the corresponding
statements for $\mathcal{L}$- and $\mathcal{R}$-classes. This
completes the proof.
\end{proof}

For $k=1,\dots,\lfloor\frac{n}{2}\rfloor$ set
$\varepsilon_k=\tau_{1,2}\tau_{3,4}\dots\tau_{2k-1,2k}$ and
let $\mathcal{H}_k$ denote the $\mathcal{H}$-class of $T$,
containing the element $\varepsilon_k$.
For $i,j\in \{3,\dots,n\}$, $i\neq j$, set
$\gamma_{i,j}=\tau_{1,2}\tau_{1,i}\tau_{1,j}\tau_{1,2}$. Note 
that, using \eqref{3} and \eqref{4}, we have
\begin{equation}\label{eq:1-2}
\tau_{1,2}\tau_{1,i}\tau_{1,j}\tau_{1,2}=
\tau_{1,2}\tau_{1,i}\tau_{i,j}\tau_{1,j}\tau_{1,2}=
\tau_{1,2}\tau_{2,j}\tau_{i,j}\tau_{2,i}\tau_{1,2}=
\tau_{1,2}\tau_{2,j}\tau_{2,i}\tau_{1,2}.
\end{equation}

\begin{lemma}\label{lm:prelude-to-Coxeter}
The elements $\gamma_{i,j}$, $i,j\in \{3,\dots,n\}$, $i\neq j$,
generate $\mathcal{H}_{1}$ as a monoid.
\end{lemma}

\begin{proof}
Let $w\in \A^+$ be such that $w\in \mathcal{H}_{1}$. Since 
$\tau_{1,2}$ is the unit element in the group $\mathcal{H}_{1}$,
we have $w=\tau_{1,2}w\tau_{1,2}$ and hence we can assume that 
$w$ has the form $\tau_{1,2}w'\tau_{1,2}$ for some $w'\in\A^+$.
We claim that $w$ is connected. Indeed, assume that $w$ is not 
connected. Then a direct calculation shows that $\varphi(w)\in\B_n$ 
has corank at least $4$. At the same time the corank 
of $\varphi(\varepsilon_1)$ is $2$. This contradicts
Lemma~\ref{lm:Green-on-Brauer}. 

We prove our lemma by induction on the length of 
$w'=\tau_{i_1,j_1}\dots\tau_{i_k,j_k}$ (note that 
$w'$ is connected since $w$ is). Because of \eqref{2} we
can always assume that $\tau_{i_s,j_s}\neq \tau_{i_{s+1},j_{s+1}}$
for all $s=1,\dots,k-1$, $\tau_{i_1,j_1}\neq \tau_{1,2}$,
and $\tau_{i_k,j_k}\neq \tau_{1,2}$. The basis of our induction will
be cases $k=0,1,2$. If $k=0,1$, then from 
\eqref{2} and \eqref{5} it follows that $w=\tau_{1,2}$, and the 
statement is obvious.

Let $k=2$. If either $1$ or $2$ occurs in both
$\{i_1,j_1\}$ and $\{i_2,j_2\}$, we are done by \eqref{eq:1-2}.
If not, without loss of generality and up to the application of
$*$ we can assume that $i_1=1$ and $i_2=2$. Then $j_1=j_2$
since $w$ is connected. Hence, using \eqref{4} we get
\begin{displaymath}
\tau_{1,2}\tau_{1,j_1}\tau_{2,j_1}\tau_{1,2}=
\tau_{1,2}\tau_{1,j_1}\tau_{1,2},
\end{displaymath}
reducing everything to the case $k=1$.

Now we proceed by induction and prove the step $k-1\Rightarrow k$,
where $k>2$. If $\{i_2,j_2\}\cap\{1,2\}\neq\varnothing$, using
\eqref{5} we can write
\begin{displaymath}
\tau_{1,2}\tau_{i_1,j_1}\tau_{i_2,j_2}\tau_{i_3,j_3}\dots
\tau_{i_k,j_k}\tau_{1,2}=
\tau_{1,2}\tau_{i_1,j_1}\tau_{i_2,j_2}\tau_{1,2}\tau_{i_2,j_2}
\tau_{i_3,j_3}\dots \tau_{i_k,j_k} \tau_{1,2}
\end{displaymath}
and the statement follows from the inductive assumption.
If $\{i_2,j_2\}\cap\{1,2\}=\varnothing$ then, using \eqref{3} if
necessary, we may assume $i_1=1$ and $j_1=j_2$. Assume first that
$j_1\in\{i_3,j_3\}$, say $j_3=j_1$. Then by \eqref{3} we have
\begin{displaymath}
\tau_{1,2}\tau_{1,j_1}\tau_{i_2,j_1}\tau_{i_3,j_1}=
\tau_{1,2}\tau_{2,i_2}\tau_{i_2,j_1}\tau_{i_3,j_1}.
\end{displaymath}
If $i_3=2$, then \eqref{4} gives
$\tau_{2,i_2}\tau_{i_2,j_1}\tau_{2,j_1}=\tau_{2,i_2}\tau_{2,j_1}$
and reduces our expression to the case $k-1$. If $i_3\neq 2$,
using \eqref{3} we have
$\tau_{2,i_2}\tau_{i_2,j_1}\tau_{i_3,j_1}=
\tau_{2,i_2}\tau_{2,i_3}\tau_{i_3,j_1}$, which reduces 
our expression to the case 
$\{i_2,j_2\}\cap\{1,2\}\neq\varnothing$, considered above.

Finally, assume that $j_1\not\in\{i_3,j_3\}$. Then, without
loss of generality we can assume $i_3=i_2$. If
$j_3=1$, then, by \eqref{4} we have
$\tau_{1,j_1}\tau_{i_2,j_1}\tau_{i_2,1}=\tau_{1,j_1}\tau_{i_2,1}$,
which reduces our expression to the case $k-1$. If $j_3\neq 1$,
using \eqref{3} we have
$\tau_{1,j_1}\tau_{i_2,j_1}\tau_{i_2,j_3}=
\tau_{1,j_1}\tau_{1,j_3}\tau_{i_2,j_3}$, which reduces 
our expression to the case 
$\{i_2,j_2\}\cap\{1,2\}\neq\varnothing$, considered above.
Now the proof is completed by induction.
\end{proof}

For $3\leq i\leq n-1$ set $\gamma_i=\gamma_{i,i+1}$. 

\begin{lemma}\label{lm:delta-is-decomposable}
Let $i,j\in\{3,\dots,n\}$, $i\neq j$.
\begin{enumerate}[(i)]
\item\label{lm:delta-is-decomposable.1}
$\gamma_{i,j}=\gamma_{j,i}$.
\item\label{lm:delta-is-decomposable.2}
$\gamma_{i,j}$ decomposes into a product of $\gamma_k$'s.
\end{enumerate}
\end{lemma}

\begin{proof}
We have:
\begin{displaymath}
\begin{array}{rcll}
\gamma_{i,j}&=& \tau_{1,2}\tau_{1,i}\tau_{1,j}\tau_{1,2}
&\\
&=& \tau_{1,2}\tau_{2,i}\tau_{1,i}\tau_{1,j}\tau_{1,2}
&\text{(by \eqref{4})}\\
&=& \tau_{1,2}\tau_{2,i}\tau_{2,j}\tau_{1,j}\tau_{1,2}
&\text{(by \eqref{3})}\\
&=& \tau_{1,2}\tau_{2,i}\tau_{2,j}\tau_{1,2}
&\text{(by \eqref{4})}\\
&=& \tau_{1,2}\tau_{1,j}\tau_{1,i}\tau_{1,2}
&\text{(by \eqref{eq:1-2})},
\end{array}
\end{displaymath}
which proves \eqref{lm:delta-is-decomposable.1}.

Because of \eqref{lm:delta-is-decomposable.1} we can assume $j>i$.
If $j=i+1$, then \eqref{lm:delta-is-decomposable.2} is obvious.
We proceed by induction on $j-i$ and assume that 
some $\gamma_{i,j}$ decomposes into a product of $\gamma_k$'s. 
We have:
\begin{equation}\label{xxx100}
\begin{array}{rcll}
\gamma_{i,j}\gamma_{j}\gamma_{i,j}&=&
\tau_{1,2}\tau_{1,i}\tau_{1,j}\tau_{1,2}\tau_{1,j}
\tau_{1,j+1}\tau_{1,2}\tau_{1,i}\tau_{1,j}\tau_{1,2}
&\text{(by \eqref{2})}\\
&=&
\tau_{1,2}\tau_{1,i}\tau_{1,j}\tau_{1,j+1}\tau_{1,2}
\tau_{1,i}\tau_{1,j}\tau_{1,2}&\text{(by \eqref{5})}\\
&=&
\tau_{1,2}\tau_{1,i}\tau_{1,j}\tau_{1,j+1}
\tau_{1,2}\tau_{2,j}\tau_{2,i}\tau_{1,2}&
\text{(by \eqref{eq:1-2})}\\
&=&
\tau_{1,2}\tau_{1,i}\tau_{1,j}\tau_{1,j+1}
\tau_{j,j+1}\tau_{2,j}\tau_{2,i}\tau_{1,2}&
\text{(by \eqref{3})}\\
&=&
\tau_{1,2}\tau_{1,i}\tau_{1,j}\tau_{j,j+1}
\tau_{2,j}\tau_{2,i}\tau_{1,2}&
\text{(by \eqref{4})}\\
&=&
\tau_{1,2}\tau_{1,i}\tau_{1,j}\tau_{j,j+1}
\tau_{i,j+1}\tau_{2,i}\tau_{1,2}&
\text{(by \eqref{3})}\\
&=&
\tau_{1,2}\tau_{1,i}\tau_{1,j}
\tau_{1,i}\tau_{i,j+1}\tau_{2,i}\tau_{1,2}&
\text{(by \eqref{3})}\\
&=&
\tau_{1,2}\tau_{1,i}\tau_{i,j+1}\tau_{2,i}\tau_{1,2}&
\text{(by \eqref{5})}\\
&=&
\tau_{1,2}\tau_{1,i}\tau_{i,j+1}\tau_{1,j+1}\tau_{1,2}&
\text{(by \eqref{3})}\\
&=&
\tau_{1,2}\tau_{1,i}\tau_{1,j+1}\tau_{1,2}&
\text{(by \eqref{4})}\\
&=&\gamma_{i,j+1}.
\end{array}
\end{equation}
The statement \eqref{lm:delta-is-decomposable.2} now follows by induction.
\end{proof}

\begin{lemma}\label{lm:Coxeter}
The elements $\gamma_i$, $i=3,\dots,n-1$, satisfy the following
relations:
\begin{enumerate}[(a)]
\item\label{cox:square} $\gamma_i^2=\tau_{1,2}$;
\item\label{cox:commute} $\gamma_i\gamma_j=\gamma_j\gamma_i$,
$\lmod i-j\rmod>1$;
\item\label{cox:braid} 
$\gamma_i\gamma_j\gamma_i=\gamma_j\gamma_i\gamma_j$,
$\lmod i-j\rmod=1$.
\end{enumerate}
\end{lemma}

\begin{proof}
We have:
\begin{displaymath}
\begin{array}{rcll}
\gamma_{i}^2&=&
\tau_{1,2}\tau_{1,i}\tau_{1,i+1}\tau_{1,2}\tau_{1,i}\tau_{1,i+1}\tau_{1,2}
&\text{(by \eqref{2})}\\
&=&
\tau_{1,2}\tau_{2,i+1}\tau_{2,i}\tau_{1,2}\tau_{1,i}\tau_{1,i+1}\tau_{1,2}
&\text{(by \eqref{eq:1-2})}\\
&=&
\tau_{1,2}\tau_{2,i+1}\tau_{2,i}\tau_{1,i}\tau_{1,i+1}\tau_{1,2}
&\text{(by \eqref{4})}\\
&=&
\tau_{1,2}\tau_{2,i+1}\tau_{1,i+1}\tau_{1,i}\tau_{1,i+1}\tau_{1,2}
&\text{(by \eqref{3})}\\
&=&
\tau_{1,2}\tau_{2,i+1}\tau_{1,i+1}\tau_{1,2}
&\text{(by \eqref{5})}\\
&=&
\tau_{1,2}\tau_{1,i+1}\tau_{1,2}
&\text{(by \eqref{4})}\\
&=&
\tau_{1,2}
&\text{(by \eqref{5})},
\end{array}
\end{displaymath}
which implies \eqref{cox:square}.

To prove \eqref{cox:commute} we may assume $j\geq i+2$. We have:
\begin{displaymath}
\begin{array}{rcll}
\gamma_i\gamma_j&=&
\tau_{1,2}\tau_{1,i}\tau_{1,i+1}\tau_{1,2}\tau_{1,j}\tau_{1,j+1}\tau_{1,2}
&\text{(by \eqref{2})}\\
&=&
\tau_{1,2}\tau_{1,i}\tau_{1,i+1}\tau_{1,2}\tau_{2,j}\tau_{1,j}\tau_{1,j+1}\tau_{1,2}
&\text{(by \eqref{4})}\\
&=&
\tau_{1,2}\tau_{1,i}\tau_{1,i+1}\tau_{i+1,j}\tau_{2,j}\tau_{1,j}\tau_{1,j+1}\tau_{1,2}
&\text{(by \eqref{3})}\\
&=&
\tau_{1,2}\tau_{1,i}\tau_{i,j}\tau_{i+1,j}\tau_{2,j}\tau_{2,j+1}\tau_{1,j+1}\tau_{1,2}
&\text{(by \eqref{3})}\\
&=&
\tau_{1,2}\tau_{1,i}\tau_{i,j}\tau_{i+1,j}\tau_{2,j}\tau_{2,j+1}\tau_{1,2}
&\text{(by \eqref{4})}\\
&=&
\tau_{1,2}\tau_{2,j}\tau_{i,j}\tau_{i+1,j}\tau_{i+1,j+1}\tau_{2,j+1}\tau_{1,2}
&\text{(by \eqref{3})}\\
&=&
\tau_{1,2}\tau_{2,j}\tau_{i,j}\tau_{i,j+1}\tau_{i+1,j+1}\tau_{2,j+1}\tau_{1,2}
&\text{(by \eqref{3})}\\
&=&
\tau_{1,2}\tau_{2,j}\tau_{2,j+1}\tau_{i,j+1}\tau_{i+1,j+1}\tau_{1,i+1}\tau_{1,2}
&\text{(by \eqref{3})}\\
&=&
\tau_{1,2}\tau_{2,j}\tau_{2,j+1}\tau_{i,j+1}\tau_{1,i}\tau_{1,i+1}\tau_{1,2}
&\text{(by \eqref{3})}\\
&=&
\tau_{1,2}\tau_{2,j}\tau_{2,j+1}\tau_{1,2}\tau_{1,i}\tau_{1,i+1}\tau_{1,2}
&\text{(by \eqref{3})}\\
&=&
\tau_{1,2}\tau_{1,j+1}\tau_{1,j}\tau_{1,2}\tau_{1,i}\tau_{1,i+1}\tau_{1,2}
&\text{(by \eqref{eq:1-2})}\\
&=&
\tau_{1,2}\tau_{1,j}\tau_{1,j+1}\tau_{1,2}\tau_{1,i}\tau_{1,i+1}\tau_{1,2}
&\text{(by Lemma~\ref{lm:delta-is-decomposable}\eqref{lm:delta-is-decomposable.1})}\\
&=&
\gamma_j\gamma_i.
&\text{(by \eqref{2})}
\end{array}
\end{displaymath}
This gives \eqref{cox:commute}. 

Finally, to prove \eqref{cox:braid} we may assume $j=i+1$. We have:
\begin{displaymath}
\begin{array}{rcll}
\gamma_{i+1}\gamma_i\gamma_{i+1}&=&
\tau_{1,2}\tau_{1,i+1}\tau_{1,i+2}\tau_{1,2}\tau_{1,i}\tau_{1,i+1}\tau_{1,2}\tau_{1,i+1}\tau_{1,i+2}\tau_{1,2}
&\text{(by \eqref{2})}\\
&=&
\tau_{1,2}\tau_{1,i+1}\tau_{1,i+2}\tau_{1,2}\tau_{1,i}\tau_{1,i+1}\tau_{1,i+2}\tau_{1,2}
&\text{(by \eqref{5})}\\
&=&
\tau_{1,2}\tau_{2,i+2}\tau_{2,i+1}\tau_{1,2}\tau_{1,i}\tau_{1,i+1}\tau_{1,i+2}\tau_{1,2}
&\text{(by \eqref{eq:1-2})}\\
&=&
\tau_{1,2}\tau_{2,i+2}\tau_{2,i+1}\tau_{i,i+1}\tau_{1,i}\tau_{1,i+1}\tau_{1,i+2}\tau_{1,2}
&\text{(by \eqref{3})}\\
&=&
\tau_{1,2}\tau_{2,i+2}\tau_{2,i+1}\tau_{i,i+1}\tau_{1,i+1}\tau_{1,i+2}\tau_{1,2}
&\text{(by \eqref{4})}\\
&=&
\tau_{1,2}\tau_{2,i+2}\tau_{2,i+1}\tau_{i,i+1}\tau_{i,i+2}\tau_{1,i+2}\tau_{1,2}
&\text{(by \eqref{3})}\\
&=&
\tau_{1,2}\tau_{2,i+2}\tau_{2,i+1}\tau_{2,i+2}\tau_{i,i+2}\tau_{1,i+2}\tau_{1,2}
&\text{(by \eqref{3})}\\
&=&
\tau_{1,2}\tau_{2,i+2}\tau_{i,i+2}\tau_{1,i+2}\tau_{1,2}
&\text{(by \eqref{5})}\\
&=&
\tau_{1,2}\tau_{1,i}\tau_{i,i+2}\tau_{1,i+2}\tau_{1,2}
&\text{(by \eqref{3})}\\
&=&
\tau_{1,2}\tau_{1,i}\tau_{1,i+2}\tau_{1,2}
&\text{(by \eqref{4})}\\
&=&
\gamma_{i,i+2}.
\end{array}
\end{displaymath}
Now \eqref{cox:braid} follows from \eqref{xxx100}.
This completes the proof.
\end{proof}

\begin{corollary}\label{c5}
\begin{enumerate}[(i)]
\item\label{c5.1} $\mathcal{H}_1\cong \Sym_{n-2}$.
\item\label{c5.2} Let $\pi\in T$ be such that 
$\pi\mathcal{D}\varepsilon_1$. Then the restriction of
$\varphi$ to $\mathcal{H}_{\pi}$ is injective.
\end{enumerate}
\end{corollary}

\begin{proof}
$\mathcal{H}_1$ contains $\varepsilon_1$ and hence is a group.
By Lemma~\ref{lm:delta-is-decomposable} and
Lemma~\ref{lm:prelude-to-Coxeter}, $\mathcal{H}_1$
is generated by $\gamma_i$, $i=3,\dots,n-1$. By
Lemma~\ref{lm:Coxeter}, $\gamma_i$'s satisfy Coxeter 
relations of type $A_{n-3}$. Hence 
$\mathcal{H}_1$ is a quotient of $\Sym_{n-2}$. However,
$\varphi(\mathcal{H}_1)$ is a maximal subgroup of
$\B_n$, which is isomorphic to $\Sym_{n-2}$ by
\cite[Theorem~1]{Maz1}. The statement \eqref{c5.1} follows.

\eqref{c5.1} implies that the restriction of 
$\varphi$ to $\mathcal{H}_{1}$ is injective. Then for
arbitrary $\pi\in T$ such that  $\pi\mathcal{D}\varepsilon_1$
the statement \eqref{c5.2} follows from Green's Lemma.
\end{proof}

To prove Theorem~\ref{th:main} we have to generalize 
the statement of Corollary~\ref{c5}\eqref{c5.2} to 
all other $\mathcal{H}$-classes. For this we will use
the following statement:

\begin{proposition}\label{pr:card-of-k}
$|\mathcal{H}_k|=(n-2k)!$ for all $k$, $1\leq k\leq\lfloor\frac{n}{2}\rfloor$.
\end{proposition}

\begin{proof}
We induct on $k$. The case $k=1$ follows from
Corollary~\ref{c5}\eqref{c5.1}. Let us prove the induction step
$k-1\Rightarrow k$. From the inductive assumption and Green's Lemma
it follows that every $\mathcal{H}$-class, which is
$\mathcal{D}$-equivalent to $\mathcal{H}_{k-1}$, has cardinality
$(n-2(k-1))!$. 

For $i=1,2,\dots,k$ set
\begin{displaymath}
\theta_i=\tau_{1,2}\tau_{3,4}\dots\tau_{2i-3,2i-2}
\tau_{2i+1,2i+2}\dots\tau_{2k-1,2k}.
\end{displaymath}
If $\pi\in \mathcal{L}_{\theta_i}$, then, using
Proposition~\ref{p1}, one shows that
$\pi\tau_{2i-1,2i}\in \mathcal{L}_{\varepsilon_k}$. 
Let $f_i:\mathcal{L}_{\theta_i}\to \mathcal{L}_{\varepsilon_k}$
denote the map $f_i(\pi)=\pi\tau_{2i-1,2i}$. This induces the map
$$f:\coprod_{i}\mathcal{L}_{\theta_i}\to \mathcal{L}_{\varepsilon_k}$$
such that the restriction of $f$ to $\mathcal{L}_{\theta_i}$
coincides with $f_i$. By Proposition~\ref{p1} we have that
$f$ is surjective. Our aim is to prove that even $f_1$ is surjective.

Let $i,j\in\{1,2,\dots,k\}$, $i\neq j$. Define
$\alpha: \mathcal{L}_{\theta_i}\to \mathcal{L}_{\theta_j}$ via
$\alpha(\pi)=\pi\tau_{2i,2j}\tau_{2i-1,2i}$; and 
$\beta: \mathcal{L}_{\theta_j}\to \mathcal{L}_{\theta_i}$ via
$\beta(\pi)=\pi\tau_{2i,2j}\tau_{2j-1,2j}$. Consider the diagram
\begin{equation}\label{diagram}
\xymatrix{
\mathcal{L}_{\theta_i}\ar[rrd]_{f_i}\ar@/_/[rrrr]_{\alpha} && && \mathcal{L}_{\theta_j}\ar[lld]^{f_j}\ar@/_/[llll]_{\beta}\\
 && \mathcal{L}_{\varepsilon_k}&&\\
}.
\end{equation}
Every element $\pi\in \mathcal{L}_{\theta_i}$ satisfies
$\pi\tau_{2j-1,2j}=\pi$ by the definition of $\theta_i$.
Every element $\pi\in \mathcal{L}_{\theta_j}$ satisfies
$\pi\tau_{2i-1,2i}=\pi$ by the definition of $\theta_j$.
Further
\begin{displaymath}
\begin{array}{rcll}
\tau_{2j-1,2j}\tau_{2i,2j}\tau_{2i-1,2i}\tau_{2j-1,2j} &=&
\tau_{2j-1,2j}\tau_{2i-1,2j-1}\tau_{2i-1,2i}\tau_{2j-1,2j} &
\text{(by \eqref{3})}\\
&=&
\tau_{2j-1,2j}\tau_{2i-1,2j-1}\tau_{2j-1,2j}\tau_{2i-1,2i} &
\text{(by \eqref{7})}\\
&=&
\tau_{2j-1,2j}\tau_{2i-1,2i} &
\text{(by \eqref{5}).}\\
\end{array}
\end{displaymath}
This implies for all $\pi\in\GL_{\theta_i}$ the following:
\begin{displaymath}
\begin{array}{rcll}
(f_j\alpha)(\pi) &=&
(f_j\alpha)(\pi\tau_{2j-1,2j})\\
&=&
\pi\tau_{2j-1,2j}\tau_{2i,2j}\tau_{2i-1,2i}\tau_{2j-1,2j}\\
&=&
\pi\tau_{2j-1,2j}\tau_{2i-1,2i}\\
&=&
\pi\tau_{2i-1,2i}\\
&=&
f_i(\pi).
\end{array}
\end{displaymath}
Hence $f_j \alpha=f_i$. Analogously one shows that $f_i\beta=f_j$.
Thus the diagram \eqref{diagram} is commutative, which implies that
the map $f_1$ is surjective. 

\begin{lemma}\label{aux1}
For any $\pi\in\GL_{\theta_1}$ there exists 
$\omega\in T$ such that $\omega\mathcal{H}\pi$, $\omega\neq \pi$,
and $f_1(\pi)=f_1(\omega)$.
\end{lemma}

\begin{proof}
Set $\omega=\pi\tau_{3,4}\tau_{1,3}\tau_{2,3}\tau_{3,4}$.
A direct calculation shows that $\varphi(\pi)\mathcal{L}\varphi(\omega)$.
Hence $\omega\in\GL_{\theta_1}$ by Corollary~\ref{cgreen}\eqref{cgreen.3}.
Further, we have $(\tau_{3,4}\tau_{1,3}\tau_{2,3}\tau_{3,4})^2=\tau_{3,4}$
by the statement analogous to that of Lemma~\ref{lm:Coxeter}\eqref{cox:square}, which implies
$\omega\mathcal{R}\pi$, that is $\omega\mathcal{H}\pi$. A direct calculation
shows that  $\varphi(\pi)\neq
\varphi(\pi)\varphi(\tau_{3,4}\tau_{1,3}\tau_{2,3}\tau_{3,4})$, and hence
$\pi\neq \omega$. On the other hand,
\begin{displaymath}
\begin{array}{rcll}
f_1(\omega)&=&
\pi\tau_{3,4}\tau_{1,3}\tau_{2,3}\tau_{3,4}\tau_{1,2}\\
&=&
\pi\tau_{3,4}\tau_{1,3}\tau_{2,3}\tau_{1,2}\tau_{3,4}
&\text{(by \eqref{7})}\\
&=&
\pi\tau_{3,4}\tau_{1,3}\tau_{1,2}\tau_{3,4}
&\text{(by \eqref{4})}\\
&=&
\pi\tau_{3,4}\tau_{1,3}\tau_{3,4}\tau_{1,2}
&\text{(by \eqref{7})}\\
&=&
\pi\tau_{3,4}\tau_{1,2}
&\text{(by \eqref{5})}\\
&=&
\pi\tau_{1,2}&\text{(since $\pi\in\GL_{\theta_1}$)}\\
&=&
f_1(\pi).
\end{array}
\end{displaymath}
\end{proof}

\begin{lemma}\label{aux2}
Assume that $\pi,\tau\in\GL_{\theta_1}$ are such that 
$f_1(\pi)\GH f_1(\tau)$. Then there exists $\eta\in\GH_{\pi}$ 
such that $f_1(\eta)=f_1(\tau)$.
\end{lemma}

\begin{proof}
If $\tau\in \GH_{\pi}$, we have nothing to prove, hence
we assume that $\tau\not\in \GH_{\pi}$. We have
$\pi\tau_{1,2}\GH\tau\tau_{1,2}$. In particular,
$\pi\tau_{1,2}\GR\tau\tau_{1,2}$. Moreover, we also have that
the corank of $\varphi(\tau\tau_{1,2})=2k$.
Then, applying $*$ to the statement of Proposition~\ref{p1}, we 
obtain that there exist 
$w,w'\in\A^+$, pairwise distinct $i_1,j_1,\dots,i_{k},j_{k}$,
and $a,b\in \{1,2,\dots,k\}$,
such that $\pi\tau_{1,2}=\tau_{i_1,j_1}\dots\tau_{i_{k},j_{k}}w$,
$\tau\tau_{1,2}=\tau_{i_1,j_1}\dots\tau_{i_{k},j_{k}}w'$,
the word $\tau_{i_{a},j_{a}}w$ is connected, and the word 
$\tau_{i_b,j_b}w'$ is connected.
Since both $\corank(\varphi(\pi))=\corank(\varphi(\tau))=2k-2$
and $\tau\not\in \GR_{\pi}$, without loss of generality we may assume
$\tau_{i_l,j_l}\pi=\pi$ for all $l=1,\dots,k-1$ and
$\tau_{i_l,j_l}\tau=\tau$ for all $l=2,\dots,k$. Then, applying
Proposition~\ref{p1}, we get some $v\in \A^+$ and
$c\in \{2,3,\dots,k\}$  such that
$\tau=\tau_{i_2,j_2}\dots\tau_{i_{k},j_{k}}v$ and
the word $\tau_{i_{c},j_{c}}v$ is connected.
Put $\eta=\tau_{i_1,j_1}\tau_{i_1,i_k}\tau$. Since
$\tau_{i_k,j_k}\tau=\tau$ by above, and
\begin{displaymath}
\tau_{i_k,j_k}\tau_{i_1,i_k}\tau_{i_1,j_1}\tau_{i_1,i_k}\tau_{i_k,j_k}=
\tau_{i_k,j_k}
\end{displaymath}
(by two applications of \eqref{5}), we have
$\eta\mathcal{L}\tau$. 

Further, since $\tau_{i_{c},j_{c}}v$ is connected, we have
$(\tau_{i_{c},j_{c}}v)(\tau_{i_{c},j_{c}}v)^*=\tau_{i_{c},j_{c}}$ by
Lemma~\ref{l3}. Using \eqref{7}, this implies $\eta\mathcal{R}
\tau_{i_1,j_1}\tau_{i_1,i_k}\tau_{i_2,j_2}\dots\tau_{i_{k},j_{k}}$.
Using \eqref{7}, we further have
\begin{displaymath}
\tau_{i_1,j_1}\tau_{i_1,i_k}\tau_{i_2,j_2}\dots\tau_{i_{k},j_{k}}=
\tau_{i_2,j_2}\dots\tau_{i_{k-1},j_{k-1}}
\tau_{i_1,j_1}\tau_{i_1,i_k}\tau_{i_{k},j_{k}}.
\end{displaymath}
Since $\tau_{i_1,j_1}\tau_{i_1,i_k}\tau_{i_{k},j_{k}}$ is connected,
by the same argument as above we have 
$\tau_{i_1,j_1}\tau_{i_1,i_k}\tau_{i_2,j_2}\dots\tau_{i_{k},j_{k}}
\mathcal{R}\tau_{i_1,j_1}\tau_{i_2,j_2}\dots\tau_{i_{k-1},j_{k-1}}$.
Hence 
$\eta\mathcal{R}\tau_{i_1,j_1}\dots\tau_{i_{k-1},j_{k-1}}$.
It follows that $\eta\mathcal{R}\pi$ and hence $\eta\mathcal{H}\pi$.

The statement now follows from the following computation
(using \eqref{5}):
\begin{multline*}
f_1(\eta)=\eta\tau_{1,2}=\tau_{i_1,j_1}\tau_{i_1,i_k}\tau\tau_{1,2}=
\tau_{i_1,j_1}\tau_{i_1,i_k}\tau_{i_1,j_1}\dots\tau_{i_{k},j_{k}}w'=\\
\tau_{i_1,j_1}\dots\tau_{i_{k},j_{k}}w'=
\tau\tau_{1,2}=f_1(\tau).
\end{multline*}
\end{proof}

Since $f_1:\mathcal{L}_{\theta_1}\to \mathcal{L}_{\varepsilon_k}$
is surjective, Lemma~\ref{aux2} implies that the restriction of
$f_1$ to $\mathcal{H}_{\theta_1}$ is a surjection on a union of
$\mathcal{H}$-classes in $\mathcal{L}_{\varepsilon_k}$. 
By Corollary~\ref{cgreen}\eqref{cgreen.3}, the number
of $\mathcal{H}$-classes in the latter union can be computed
in the semigroup $\B_n$ via $\varphi$, and it is easy to see
that it equals $\binom{n-(2k-2)}{2}$.

We know by induction that $|\mathcal{H}_{\theta_1}|=(n-2(k-1))!$.
Hence, taking into account Lemma~\ref{aux1} and Green's Lemma, we compute:
\begin{equation}\label{for12}
\lmod\mathcal{H}_{f_1(\theta_1)}\rmod\leq
\dfrac{1}{\binom{n-(2k-2)}{2}}\cdot\frac{(n-2(k-1))!}{2}=
(n-2k)!.
\end{equation}
Since $\lmod\mathcal{H}_{\varphi(f_1(\theta_1))}\rmod=(n-2k)!$ by
\cite[Theorem~1]{Maz1}, \eqref{for12} and 
Corollary~\ref{cgreen}\eqref{cgreen.3} imply
$\lmod\mathcal{H}_{f_1(\theta_1)}\rmod=(n-2k)!$. This forces
$|\mathcal{H}_{k}|=(n-2k)!$ by Green's Lemma and the statement
follows by induction.
\end{proof}

\begin{proof}[Proof of Theorem~\ref{th:main}]
Let $1\leq k\leq \lfloor\frac{n}{2}\rfloor$. By
Proposition~\ref{pr:card-of-k} we have 
$|\mathcal{H}_k|=(n-2k)!$. By \cite[Theorem~1]{Maz1} we have
$|\varphi(\mathcal{H}_k)|=(n-2k)!$ as well. Hence the restriction
of $\varphi$ to $\mathcal{H}_k$ is injective. From Green's Lemma
it follows that the restriction of $\varphi$ to $\mathcal{H}_{\pi}$
is injective for every $\pi\in T$ such that $\pi\mathcal{D}\varepsilon_k$.
From Corollary~\ref{cgreen}\eqref{cgreen.3} it therefore follows that
$\varphi$ is injective, and hence bijective. This completes the
proof.
\end{proof}

%%%%%%%%%%%%%%%%%%%%%%%%%%%%%%%%%%%%%%%%%%%%%%%%%%%%%%%%%%%%%%%%%%%%%%%%%%%%%%%%%%%%%%%%%%%%%%%%%

\section{Combinatorial applications}\label{sec:application}

%%%%%%%%%%%%%%%%%%%%%%%%%%%%%%%%%%%%%%%%%%%%%%%%%%%%%%%%%%%%%%%%%%%%%%%%%%%%%%%%%%%%%%%%%%%%%%%%%

\subsection{Connected sequences}\label{applcat}

Two elements, $\{i,j\}$ and $\{k,l\}$, of $\binom{{\bn}}{2}$ are 
said to be \emph{connected} provided that 
$\{i,j\}\cap\{k,l\}\ne\varnothing$. A {\em connected sequence} then 
is a non-empty sequence, $\{i_1,j_1\}$, $\{i_2,j_2\}$,\dots,
$\{i_m,j_m\}$, of elements from $\binom{{\bn}}{2}$ such that 
$\{i_l,j_l\}$ and
$\{i_{l+1},j_{l+1}\}$ are connected for all $l=1,\dots,m-1$. 
Two connected sequences will be called {\em equivalent} provided
that one of them can be obtained from the other one by a finite
number of the following operations:
\begin{enumerate}[(I)]
\item\label{seq:2} replacing the fragment $\{i,j\},\{i,j\}$ by $\{i,j\}$
and vice versa;
\item\label{seq:3}
replacing the fragment $\{i,j\},\{j,k\},\{k,l\}$ by $\{i,j\},\{i,l\},\{k,l\}$ and 
vice versa, where $i\neq l$;
\item\label{seq:4}
replacing the fragment $\{i,j\},\{j,k\},\{k,i\}$ by $\{i,j\},\{k,i\}$ and 
vice versa;
\item\label{seq:5}
replacing the fragment $\{i,j\},\{j,k\},\{i,j\}$ by $\{i,j\}$ and 
vice versa.
\end{enumerate}
It is obvious that each of the operations \eqref{seq:2}-\eqref{seq:5},
applied to a connected sequence, produces a new connected sequence.
As an immediate corollary of  Theorem~\ref{th:main} we have the 
following result:

\begin{proposition}\label{pr:application-sequences}
Let $n\in\{2,3,\dots\}$.
\begin{enumerate}[(i)]
\item There exist only finitely many, namely $\dfrac{n(n-1)n!}{4}$, 
equivalence classes of connected sequences.
\item For all $\{i,j\},\{k,l\}\in \binom{{\bn}}{2}$ the number of connected
sequence, whose first element is $\{i,j\}$ and whose last element is
$\{k,l\}$, equals $(n-2)!$.
\end{enumerate}
\end{proposition}

\begin{proof}
Let $S$ denote the set of all equivalence classes of connected
sequences. Define a semigroup structure on $S\cup\{0\}$ as follows:
$0$ is the zero element of $S\cup\{0\}$, and for $f,g\in S$
\begin{displaymath}
f\cdot g=\begin{cases}
fg, & \text{$fg$ is connected}\\
0, & \text{otherwise}.
\end{cases}
\end{displaymath}

Let $\overline{\B}$ denote the Rees quotient of $\B_n\setminus\Sym_n$ 
modulo  the ideal, containing all elements of corank at least $4$. 
By Theorem~\ref{th:main}, mapping $\sigma_{i,j}$ to the connected
sequence $\{i,j\}$ defines an epimorphism, $\psi$, from
$\overline{\B}$ to $S\cup\{0\}$. On the other hand, from the
definition of the equivalence relation on the connected sequences
we have that, mapping $\{i,j\}$ to $\sigma_{i,j}$ defines an
epimorphism, $\psi':S\cup\{0\}\to \overline{\B}$. Thus 
$\psi$ and $\psi'$ induce  a pair of mutually inverse bijections
between the set of all elements in $\B_n$ of corank $2$ and the
set of equivalence classes of connected sequences. The claim
now follows by a direct computation in $\B_n$.
\end{proof}

It might be interesting to find a purely combinatorial
proof for the statement of Proposition~\ref{pr:application-sequences}.

%%%%%%%%%%%%%%%%%%%%%%%%%%%%%%%%%%%%%%%%%%%%%%%%%%%%%%%%%%%%%%%%%%%%%%%%%%%%%%%%%%%%%%%%%%%%%

\subsection{Paths in the graph $\Gamma_n$}

\begin{figure}
\begin{displaymath}
\xymatrix{
(1,2)\ar@{-}[r]\ar@{-}[d]\ar@{-}[dr]\ar@/^1pc/@{-}[rr] & 
(1,3)\ar@{-}[r]\ar@{-}[dr]\ar@{-}[dl] & 
(1,4)\ar@{-}[d]\ar@{-}[dl]\\
(2,3)\ar@{-}[r]\ar@/_1pc/@{-}[rr] & (2,4)\ar@{-}[r] & (3,4)\\
}
\end{displaymath}
\caption{The graph $\Gamma_4$.}\label{fig:appl-graphs}
\end{figure}
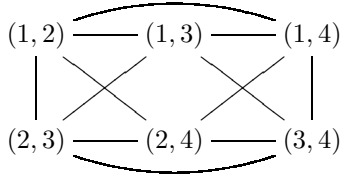

There is another interesting combinatorial interpretation of
the elements of $\B_n$ of corank $2$. Consider a non-oriented graph, 
$\Gamma_n$, whose vertex set is $\binom{\bn}{2}$, and such that 
two vertices, $\{i,j\}$ and $\{k,l\}$ are connected by an edge if
and only if $\{i,j\}\cap\{k,l\}\ne\varnothing$. The graph $\Gamma_4$
is shown on Figure~\ref{fig:appl-graphs}.

Obviously, the paths in $\Gamma_n$ can be interpreted as
connected sequences as defined in the previous subsection.
Then the equivalence relation on the connected sequences, defined
by the operations \eqref{seq:2}-\eqref{seq:5}, has the following
interpretation in terms of the graph $\Gamma_n$:
\begin{enumerate}[(I)]
\item the trivial path in vertex is an idempotent;
\item if the full subgraph of $\Gamma_n$, corresponding to 
a quadruple of vertices  has the form
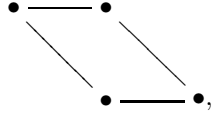
\begin{displaymath}
\xymatrix{
\bullet\ar@{-}[r]\ar@{-}[rd] & \bullet\ar@{-}[rd] & \\
& \bullet\ar@{-}[r] & \bullet, \\
}
\end{displaymath}
then the paths of length $2$ in this subgraph with the same initial
and the same terminal points are equivalent;
\item for any triple $\{i,j\}$, $\{j,l\}$, $\{i,l\}$ of vertices
the paths in the full subgraph $\Gamma_n$, corresponding to these
vertices, with the same initial and the same terminal points 
are equivalent;
\item the path, consisting of going along the same edge in
two different directions coincides with the trivial path in the
starting point.
\end{enumerate}
These relations generate an equivalence relation on the set of
all paths in $\Gamma_n$. From Proposition~\ref{pr:application-sequences}
it thus follows that the number of non-equivalent paths in
$\Gamma_n$ equals $\frac{n(n-1)n!}{4}$, and the number 
of non-equivalent loops in each point equals $(n-2)!$.

%%%%%%%%%%%%%%%%%%%%%%%%%%%%%%%%%%%%%%%%%%%%%%%%%%%%%%%%%%%%%%%%%%%%%%%%%%%%%%%%%%%%%%%%%%%%%%%%%%%%%%%%%%

\subsection{The maximal length of an element from $\B_n\setminus\Sym_n$}\label{indeed-nice}

For $w\in\B_n\setminus\Sym_n$ we define the {\em length} 
$\ls(w)$ of $w$ as the length of the shortest possible
presentation of $w$ as a product of the generators $\sigma_{i,j}$'s.
For $w\in \A_n^+$ we define the {\em length} 
$\len(w)$ of $w$ as the length of
presentation of $w$ as a product of the generators $\tau_{i,j}$'s.
The aim of this subsection is to prove the following statement about
the maximal value $f(n)$ of $\ls(w)$ on $w\in \B_n\setminus\Sym_n$.

\begin{theorem}\label{th:indeed-nice}
Let $n\geq 2$. Then $f(n)=\lfloor\frac{3n}{2}\rfloor-2$.
\end{theorem}

For the proof of Theorem~\ref{th:indeed-nice} we will need several
auxiliary statements. Set $g(n)=\lfloor\frac{3n}{2}\rfloor-2$.
We will show that $f(n)\leq g(n)$. For $n\neq 3$ we will then
find an element,  $w\in\GH_{\sigma_{1,2}}$, such that 
$\ls(w)=g(n)$. For $n=3$ we have $\ls(\sigma_{1,2}\sigma_{2,3})=2=g(3)$.
By Theorem~\ref{th:main} we have $\B_n\setminus\Sym_n\simeq T$ and hence
in the sequel we can work with the semigroup $T$ and the generators
$\tau_{i,j}$'s. The function $\ls$ on $T$ defined in the obvious way,
and we consider all elements of $\A^+$ as elements of $T$ via the
natural projection.

\begin{lemma}\label{lm:1-but-not-2}
Let $u_i\in \{3,\dots,n\}$, $1\leq i\leq k$. Then
\begin{equation*}
\tau_{1,2}\tau_{1,u_1}\dots\tau_{1,u_k}\tau_{1,2}=\tau_{1,2}\tau_{2,u_1}\dots\tau_{2,u_k}\tau_{1,2}.
\end{equation*}
\end{lemma}

\begin{proof}
Because of \eqref{2} we may assume that $u_i\neq u_{i+1}$ for
all $i=1,\dots,k-1$. We have:
\begin{displaymath}
\begin{array}{rcll}
\tau_{1,2}\tau_{1,u_1}\dots\tau_{1,u_k}\tau_{1,2} &=&
\tau_{1,2}\tau_{2,u_1}\tau_{1,u_1}\tau_{1,u_2}\dots\tau_{1,u_k}\tau_{1,2} &
\text{(by \eqref{4})}\\
&=&
\tau_{1,2}\tau_{2,u_1}\tau_{2,u_2}\tau_{1,u_2}\dots\tau_{1,u_k}\tau_{1,2} &
\text{(by \eqref{3})}\\
&\dots&\\
&=&
\tau_{1,2}\tau_{2,u_1}\tau_{2,u_2}\dots\tau_{2,u_k}\tau_{1,u_k}\tau_{1,2} &
\text{(by \eqref{3})}\\
&=&
\tau_{1,2}\tau_{2,u_1}\tau_{2,u_2}\dots\tau_{2,u_k}\tau_{1,2} &
\text{(by \eqref{4}).}\\
\end{array}
\end{displaymath}
\end{proof}

\begin{corollary}\label{cor:1-but-not-2}
Let $u_i\in \{2,\dots,n\}$, $1\leq i\leq k$. Then there exist elements 
$v_i\in\{1,\dots,n\}\setminus\{2\}$, $1\leq i\leq k$, such that
\begin{equation*}
w=\tau_{1,2}\tau_{1,u_1}\dots\tau_{1,u_k}\tau_{1,2}=
\tau_{1,2}\tau_{2,v_1}\dots\tau_{2,v_k}\tau_{1,2}.
\end{equation*}
\end{corollary}

\begin{proof}
If $2\notin\{u_1,\dots,u_k\}$ then the statement follows from Lemma~\ref{lm:1-but-not-2}. Otherwise let
$u_{i_1}=u_{i_2}=\dots=u_{i_p}=2$, $i_1<i_2<\dots<i_p$, be all 
occurrences of $2$  among $u_1,\dots,u_k$. Set
$u_{0}=u_{k+1}=2$, $i_0=0$, $i_{p+1}=k+1$. The statement now follows 
by applying Lemma~\ref{lm:1-but-not-2} to each element
$\tau_{1,u_{i_j}}\tau_{1,u_{i_j+1}}\dots\tau_{1,u_{i_{j+1}}}$,
$j=0,\dots,p$.
\end{proof}

\begin{lemma}\label{lm:let-us-write-1}
Let $w\in \A^+$ be such that $w\in \GH_1$. Assume that 
$\len(w)=\ls(w)\geq 4$ and set $m=\len(w)-2$. Then there exist
$u_i\in\{2,\dots,n\}$, $i=1,\dots,m$, such that
$w=\tau_{1,2}\tau_{1,u_1}\dots\tau_{1,u_m}\tau_{1,2}$.
\end{lemma}

\begin{proof}
We induct on $m=\len(w)-2$. If $m=2$, we have
$w=\tau_{1,2}\tau_{a,b}\tau_{c,d}\tau_{1,2}$. Without loss of
generality we may assume $a,d\in\{1,2\}$. If $a=d=1$, we have nothing
to prove. If $a=d=2$, the statement follows from 
Lemma~\ref{lm:1-but-not-2}. If $a\neq d$, using \eqref{4} we see 
that $\ls(w)<4$, a contradiction. 

Now we prove the induction step $m-1\Rightarrow m$.
Let $w=\tau_{1,2}\tau_{i_1,j_1}\dots\tau_{i_m,j_m}\tau_{1,2}$. 
Without loss of generality we may assume $i_1\in\{1,2\}$. 

{\bf Case~1:} $i_1=1$. If $1\in\{i_l,j_l\}$ for all  $l\leq m$,
we have nothing to prove. Otherwise let $p<m$ be such that
$1\in\{i_l,j_l\}$ for all  $l\leq p$ and 
$1\not\in \{i_{p+1},j_{p+1}\}$.  Thus, without loss of
generality we may assume $i_l=1$ for all  $l\leq p$.
If $j_p=2$, the statement follows from the inductive assumption.

Assume that $j_p\neq 2$. Without loss of
generality we may write
\begin{displaymath}
w=\tau_{1,2}\tau_{1,j_1}\dots\tau_{1,j_p}
\tau_{j_p,j_{p+1}}\dots\tau_{j_p,j_{p+q}}
\tau_{j_{p+q},j_{p+q+1}}\dots\tau_{i_m,j_m}\tau_{1,2}.
\end{displaymath}
Observe that $j_{p+q+1}\neq j_p$ since $\len(w)=\ls(w)$.
 
Now we apply consequently the relation \eqref{3} starting from 
$\tau_{j_p,j_{p+q}}$ and moving to the left until we reach the
element $\tau_{j_p,j_{p+2}}$. We will get
\begin{displaymath}
w=\tau_{1,2}\tau_{1,j_1}\dots\tau_{1,j_p}
\tau_{j_p,j_{p+1}}\tau_{j_{p+1},j_{p+q+1}}\dots\tau_{j_{p+q-1},j_{p+q+1}}
\tau_{j_{p+q},j_{p+q+1}}\dots\tau_{i_m,j_m}\tau_{1,2}.
\end{displaymath}
If $j_{p+q+1}=1$, then we can reduce the length of $w$ by
\eqref{4}, a contradiction. Otherwise, using \eqref{3} we
can change $\tau_{j_p,j_{p+1}}$ into $\tau_{1,j_{p+q+1}}$. 
Since we have not changed the length of $w$, the proof in Case~1 
is now completed by induction on $p$.

{\bf Case~2:} $i_1=2$. Analogously to Case~1, we get the existence of
$u_i\in \{1,3,\dots,n\}$, $i=1,\dots,m$, such that 
$w=\tau_{1,2}\tau_{2,u_1}\dots\tau_{2,u_m}\tau_{1,2}$.
Now the claim follows from Corollary~\ref{cor:1-but-not-2}.
\end{proof}

Note that for every $w\in{\GH_1}$ there exists a unique 
permutation, $\pi\in\Sym_{n-2}$, such that 
\begin{equation}\label{perm}
 w=\bigl\{\{1,2\},\{1',2'\},\{k,\pi(k-2)'+2\}_{k\ne 1,2}\bigr\}.
\end{equation}
Let $({i'}_1^{(1)},\dots,{i'}_{p_1}^{(1)})$,\dots,
$({i'}_1^{(s)},\dots,{i'}_{p_s}^{(s)})$ be a complete list of cycles
of $\pi$, which have length at least $2$. Set $i_a^{(b)}={i'}_a^{(b)}+2$
for all possible $a,b$. Then, by \eqref{cycleformula}, we have the 
following decomposition of $w$:
\begin{equation*}
w=\tau_{1,2}\tau_{1,i_{1}^{(1)}}\dots\tau_{1,i_{p_1}^{(1)}}\tau_{1,2}
\dots\tau_{1,2}\tau_{1,i_{1}^{(s)}}\dots\tau_{1,i_{p_s}^{(s)}}\tau_{1,2}.
\end{equation*}
We will call this decomposition a {\em cyclic} decomposition of $w$.
We will also say that $\tau_{1,2}$ is the cyclic decomposition of
$\tau_{1,2}$. In the obvious way we now define {\em cycles} in 
${\GH_1}$.

\begin{lemma}\label{cycle}
Let $w\in {\GH_1}$ be a non-trivial cycle. Then $\ls(w)$ equals the 
length of the cyclic decomposition of $w$.
\end{lemma}

\begin{proof}
By Lemma~\ref{lm:let-us-write-1}, there exist
$u_1,\dots,u_{\ls(w)-2}\in\{2,\dots,n\}$ such that 
\begin{displaymath}
w=\tau_{1,2}\tau_{1,u_1}\dots\tau_{1,u_{\ls(w)-2}}\tau_{1,2}.
\end{displaymath}
But then all elements from $\{3,\dots,n\}$, moved by the cycle
$w$ should obviously occur among  $u_1,\dots,u_{\ls(w)-2}$. 
The claim now follows from the formula \eqref{cycleformula}.
\end{proof}

Finally, for Theorem~\ref{th:main} we obtain that for pairwise
distinct $1,u,u_1,\dots,u_k$ and for any $l\in{\bf k}$ holds
\begin{equation}\label{eq:useful}
\tau_{1,u}\tau_{1,u_1}\dots\tau_{1,u_k}\tau_{1,u}=
\tau_{1,u}\tau_{1,u_l}\dots\tau_{1,u_k}\tau_{1,u_{1}}\dots
\tau_{1,u_{l-1}}\tau_{1,u}.
\end{equation}

\begin{lemma}\label{lm:on-commutativeness}
Let $1,a,u,u_1,\dots,u_k$ be pairwise distinct. Then
\begin{equation*}
\tau_{1,a}\tau_{1,u}\tau_{1,u_1}\dots\tau_{1,u_k}\tau_{1,u}=
\tau_{1,a}\tau_{1,u_1}\dots\tau_{1,u_k}\tau_{1,a}\tau_{1,u}.
\end{equation*}
\end{lemma}

\begin{proof}
\begin{displaymath}
\begin{array}{lcl}
\tau_{1,a}\tau_{1,u}\tau_{1,u_1}\dots\tau_{1,u_k}\tau_{1,u} 
&=& \text{(by \eqref{5})}\\
\tau_{1,a}\tau_{1,u}\tau_{1,u_1}\dots\tau_{1,u_k}\tau_{1,a}
\tau_{1,u_k}\tau_{1,u} &=&\text{(by \eqref{eq:useful})}\\
\tau_{1,a}\tau_{1,u_1}\dots\tau_{1,u_k}\tau_{1,u}\tau_{1,a}
\tau_{1,u_k}\tau_{1,u} &=&\text{(by \eqref{eq:useful})}\\
\tau_{1,a}\tau_{1,u_1}\dots\tau_{1,u_k}\tau_{1,u}\tau_{1,u_k}
\tau_{1,a}\tau_{1,u} &=&\text{(by \eqref{5})}\\
\tau_{1,a}\tau_{1,u_1}\dots\tau_{1,u_k}\tau_{1,a}\tau_{1,u}. &
\end{array}
\end{displaymath}
\end{proof}

\begin{lemma}\label{lm:cyclic-decomposition}
Let $w\in {\GH_1}$. Then the cyclic decomposition of $w$ is of length $\ls(w)$.
\end{lemma}

\begin{proof}
We induct on $\ls(w)$. If $\ls(w)\leq 3$ then the statement is trivial
since, by \eqref{2} and \eqref{5}, the only possibility is $w=\tau_{1,2}$. 
Let us now prove the induction step $m+1\Rightarrow m+2$.

Let $w\in {\GH_1}$ be such that $\ls(w)=m+2$. By 
Lemma~\ref{lm:let-us-write-1}, we may write
$w=\tau_{1,2}\tau_{1,u_1}\dots\tau_{1,u_m}\tau_{1,2}$ for some
$u_i\in\{2,\dots,n\}$, $i=1,\dots,m$.  We set $u_0=u_{m+1}=2$.  
If all of $u_i$'s are pairwise 
distinct, the word $w$ is a cycle and the statement follows. 
Suppose now that there are some repetitions among $u_0,u_1,\dots,u_m$. 
Take the leftmost element, which repeats in this series, say $u_i=u$. 
Let $j>i$ be minimal possible such that $u_j=u$. Consider the element
\begin{displaymath}
w'=\tau_{1,u_i}\dots\tau_{1,u_j}=\tau_{1,u}
\tau_{1,u_{i+1}}\dots\tau_{1,u_{j-1}}\tau_{1,u}\in 
\GH_{\tau_{1,u}}.
\end{displaymath}
Since $\ls(w)=m+2$, $\ls(w')=j-i+1<m+2$. Hence, using the inductive
assumption, the cyclic decomposition of $w'$ has length $j-i+1$.
Without loss of generality we hence may assume that the subword
$w'$ of $w$ already coincides with the corresponding cyclic
decomposition, that is is a cycle.

Now we claim that $u_{i-1},u_{i},\dots,u_{j-1}$ are pairwise 
distinct. Indeed, if not, then $u_{i-1}$
coincides with one of $u_{i+1},\dots,u_{j-1}$.
Then, applying \eqref{eq:useful} we can obtain the fragment
$\tau_{1,u_{i-1}}\tau_{1,u_i}\tau_{1,u_{i-1}}$,
which can be shortened by \eqref{5}, a contradiction.
Hence, Lemma~\ref{lm:on-commutativeness} gives
\begin{multline*}
\tau_{1,u_{i-1}}\tau_{1,u_{i}}\tau_{1,u_{i+1}}
\dots\tau_{1,u_{j-2}}\tau_{1,u_{j-1}}
\tau_{1,u_{i}}\tau_{1,u_{j+1}}=\\=
\tau_{1,u_{i-1}}\tau_{1,u_{i+1}}\tau_{1,u_{i+2}}
\dots \tau_{1,u_{j-1}}\tau_{1,u_{i-1}}
\tau_{1,u_{i}}\tau_{1,u_{j+1}}.
\end{multline*}
This operation makes the index of the first letter with repetition smaller.
Hence, applying this procedure as many times as necessary, we may assume
that $i=0$. This means that $w$ is a product of a cycle with some
element $v$ from ${\GH_1}$ of strictly smaller length. By inductive 
assumption, we may assume that $v$ is written in its cyclic decomposition.
We are left to prove that non of the elements $u_1,\dots,u_{j-1}$
occurs among cycles in $v$. Assume that some of these elements does occur.
Then, using \eqref{eq:useful}, we may assume that this is $u_{j-1}$.
At the same time, the cycles of any cyclic decomposition commute
and hence, using this and \eqref{eq:useful} we may assume that 
$u_{j-1}=u_{j+1}$. In this case we can make $w$ shorter by applying
\eqref{5}, a contradiction. This completes the proof.
\end{proof}

\begin{proposition}\label{pr:H_1}
Let $w\in {\GH_1}$. Then $\ls(w)\leq g(n)$. Moreover, if $n\geq 4$
then there exists $v\in\GH_1$ such that $\ls(v)=g(n)$.
\end{proposition}

\begin{proof}
If $w=\tau_{1,2}$ or $n\leq 3$ then the statement is obvious. 
Suppose now that $w\ne\tau_{1,2}$ and $n\geq 4$.  Let 
$\pi\in \Sym_{n-2}$ be the permutation, which corresponds to $w$
by \eqref{perm}. Let $c$ and $s$ be the number of non-trivial and trivial
cycles in $\pi$ respectively. From Lemma~\ref{lm:cyclic-decomposition}
it follows that $\ls(w)=(n-2)-s+c+1$. 

{\bf Case~1:} $n=2k$, $k\in\mathbb{N}$. Then $\ls(w)=(2k-1)+c-s\leq (2k-1)+\frac{2k-2}{2}=3k-2=g(n)$ and the equality holds if and only if 
$\pi$ contains $k-1$ transpositions.

{\bf Case~2:} $n=2k+1$, $k\in\mathbb{N}$. If $s=0$ then there should 
exist a cycle in $\pi$ of length at least $3$. Then $\ls(w)=2k+c\leq 2k+\frac{n-2-3}{2}+1=3k-1=g(n)$ and the equality holds if and only if 
$\pi$ contains one cycle of length $3$ and $k-2$ transpositions. 
If $s\geq 1$ then 
$\ls(w)=2k+c-s\leq 2k-1+c\leq 2k-1+\frac{n-2-1}{2}=3k-2<g(n)$.
The proof is complete.
\end{proof}

\begin{lemma}\label{lm:aux}
Let $i_1,j_1,\dots,i_k,j_k$ be pairwise distinct elements from $\bn$.
Then there exists a word, $\mu\in\A^+$, such that $\len(\mu)\leq 2k$,
\begin{displaymath}
\tau_{1,2}\dots\tau_{2k-1,2k}\mu\in\GL_{\tau_{i_1,j_1}\dots\tau_{i_k,j_k}},
\end{displaymath}
and $\{m,m'\}\in\varphi(\mu)$, 
$m\in\bn\setminus\{i_1,j_1,\dots,i_k,j_k,1,2,\dots,2k-1,2k\}$.
\end{lemma}

\begin{proof}
For $a,b,c,d\in\bn$, $a\ne b$, $c\ne d$, set 
\begin{displaymath}
\mu_{a,b,c,d}=\begin{cases}
\tau_{a,c}\tau_{c,d}, & \text{if 
$\{a,b\}\cap\{c,d\}=\varnothing$}\\
\tau_{c,d}, & \text{otherwise}. 
\end{cases}
\end{displaymath}
Note that $\{m,m'\}\in\varphi(\mu_{a,b,c,d})$ provided that $m\ne a,b,c,d$.

Now a direct calculation implies that 
$\tau_{a,b}\mu_{a,b,c,d}\mu_{c,d,a,b}=\tau_{a,b}$ for all
$a,b,c,d$ such that $a\ne b$ and $c\ne d$. In particular, it follows that
for any $w\in T$ such that $w\tau_{a,b}=w$ we have that the
coranks of the elements $\varphi(w)$ and $\varphi(w\mu_{a,b,c,d})$
coincide.

In particular, the element $\varphi(\mu_1)$, where 
$\mu_1=\tau_{1,2}\dots\tau_{2k-1,2k}\mu_{1,2,i_1,j_1}$, has corank
$2k$. Further, $\mu_1$ satisfies $\mu_1 \tau_{i_1,j_1}=\mu_1$ by the 
definition of $\mu_{1,2,i_1,j_1}$. Hence there exist pairwise distinct
$a_1,b_1,\dots,a_{k-1},b_{k-1}$  from
$\mathbf{n}\setminus\{i_1,j_1\}$ such that
\begin{displaymath}
\mu_1\in \GL_{\tau_{i_1,j_1}\tau_{a_1,b_1}\dots\tau_{a_{k-1},b_{k-1}}}.
\end{displaymath}
Analogously, the element  $\varphi(\mu_2)$, where
$\mu_2=\mu_1\mu_{a_1,b_1,i_2,j_2}$, also has corank $2k$. The
element $\mu_2$ satisfies $\mu_2\tau_{i_1,j_1}=\mu_2$ and
$\mu_2\tau_{i_2,j_2}=\mu_2$ by the definition of
$\mu_{a_1,b_1,i_2,j_2}$. Hence, there exist pairwise distinct
$c_1,d_1,\dots,c_{k-2},d_{k-1}$ from
$\mathbf{n}\setminus\{i_1,j_1,i_2,j_2\}$ such that
\begin{displaymath}
\mu_2\in
\GL_{\tau_{i_1,j_1}\tau_{i_2,j_2}\tau_{c_1,d_1}\dots\tau_{c_{k-2},d_{k-2}}}.
\end{displaymath}
Continuing this process $k-2$ more steps we will construct the element 
$\mu_k$ with desired properties.
\end{proof}

\begin{proof}[Proof of Theorem~\ref{th:indeed-nice}]
Let $n\geq 4$. Then, by Proposition~\ref{pr:H_1}, there is
$w\in T$ such that $\ls(w)=g(n)$. For $n=2,3$ an example of
$w\in T$ such that $\ls(w)=g(n)$ was constructed immediately after
the formulation of Theorem~\ref{th:indeed-nice}. Hence
we are left to show that for any $w\in T$ we have
$\ls(w)\leq g(n)$. Without loss of generality it is even enough
to consider those $w$ for which $\tau_{1,2}w=w$.

Let now $w\in T$ be such that $\tau_{1,2}w=w$. Assume first
that $\corank(\varphi(w))=2$. Then there exists a unique
$\{i,j\}$, $i\neq j$, such that $w\tau_{i,j}=w$. Without loss of
generality we have one of the following cases:

{\bf Case 1:} $\{i,j\}=\{1,2\}$. Then the statement follows
from Proposition~\ref{pr:H_1}.

{\bf Case 2:} $i=1$ and $j\ne 2$. Then, applying \eqref{2} and 
\eqref{5}, we obtain $w=w\tau_{1,j}\tau_{1,2}\tau_{1,j}$. 
Set $w'=w\tau_{1,j}\tau_{1,2}=w\tau_{1,2}$ and we have
$w=w'\tau_{1,j}$ by \eqref{5}. It follows that $w'\GH\tau_{1,2}$.
Consider the cyclic decomposition of  $w'$.
Assume that the cycles, occurring in this decomposition,
do not move the element $j$.
Then, by Lemma~\ref{lm:cyclic-decomposition}
and Proposition~\ref{pr:H_1}, we have that the length of
this decomposition is at most $g(n-1)$. Since 
$w=w'\tau_{1,j}$, we have $\ls(w)\leq 1+g(n-1)\leq g(n)$.
Assume now that there is a cycle in $w'$, which moves $j$.
Using \eqref{eq:useful} and the fact that the cycles in 
the cyclic decomposition commute, we may write 
$w'=w''\tau_{1,j}\tau_{1,2}$, for some $w''$ such that
$\ls(w'')=\ls(w')-2$. Since $w=w'\tau_{1,j}$, we have
$w=w''\tau_{1,j}$ by \eqref{5}, and thus
$\ls(w)<\ls(w')\leq g(n)$.

{\bf Case 3:} $\{i,j\}\cap\{1,2\}=\varnothing$. Then we can write
$w=w'\tau_{1,i}\tau_{i,j}$, $w'\in \GH_{\tau_{1,2}}$.
Consider again the cyclic decomposition of $w'$. If both $i$ and $j$
are not moved by all the cycles, we have
$\ls(w)\leq 2+g(n-2)\leq g(n)$. If $i$ is moved, then, as in the 
Case~2, we can write $w'=w''\tau_{1,i}\tau_{1,2}$ for some
$w''$ such that $\ls(w'')=\ls(w')-2$ and, using \eqref{5}, we obtain
$\ls(w)\leq \ls(w'')+2=\ls(w')\leq g(n)$. Finally, let us assume that
$i$ is not moved but $j$ is moved. Assume that 
$x_1,\dots,x_p,j$, $p>0$, is a cycle  in the cyclic decomposition of
$w'$. Then, using \eqref{eq:useful} and the fact that the cycles in 
the cyclic decomposition commute, we may assume that this cyclic
decomposition has the following form:
\begin{displaymath}
\tau_{1,2}\dots\tau_{1,2}\tau_{1,x_1}\dots\tau_{1,x_p}\tau_{1,j}\tau_{1,2}.
\end{displaymath}
Now we can compute the following expression, containing 
the last cycle of this decomposition:
\begin{displaymath}
\begin{array}{lcl}
\tau_{1,2}\tau_{1,x_1}\dots\tau_{1,x_p}\tau_{1,j}\tau_{1,2}\tau_{1,i}\tau_{i,j}
&=& \text{(by \eqref{3})}\\
\tau_{1,2}\tau_{1,x_1}\dots\tau_{1,x_p}\tau_{1,j}\tau_{1,2}\tau_{2,j}\tau_{i,j}
&=& \text{(by \eqref{4})}\\
\tau_{1,2}\tau_{1,x_1}\dots\tau_{1,x_p}\tau_{1,j}\tau_{2,j}\tau_{i,j} 
&=& \text{(by \eqref{3})}\\
\tau_{1,2}\tau_{1,x_1}\dots\tau_{1,x_p}\tau_{2,x_p}\tau_{2,j}\tau_{i,j} 
&=& \text{(by \eqref{3})}\\
&\dots&\\
\tau_{1,2}\tau_{1,x_1}\tau_{2,x_1}\dots\tau_{2,x_{p-1}}\tau_{2,x_p}\tau_{2,j}\tau_{i,j}
&=& \text{(by \eqref{4})}\\
\tau_{1,2}\tau_{2,x_1}\dots\tau_{2,x_{p-1}}\tau_{2,x_p}\tau_{2,j}\tau_{i,j}.
\end{array}
\end{displaymath}
It follows that $\ls(w)\leq \ls(w')\leq g(n)$.

Assume now that $\corank(\varphi(w))=2k$, $k>1$. We may further assume that
$\tau_{1,2}\dots\tau_{2k-1,2k}w=w$ and 
$w\tau_{i_1,j_1}\dots\tau_{i_k,j_k}=w$ for some
pairwise distinct elements $i_1,j_1,\dots,i_k,j_k$. Without loss of 
generality we may assume that we have one of  the following
cases:

{\bf Case~1:} $\{i_1,j_1,\dots,i_k,j_k\}=\{1,\dots,2k\}$.
Since the map 
\begin{displaymath}
\begin{array}{rcl}
\GH_{\tau_{1,2}\dots\tau_{2k-1,2k}} & \to & 
\GH_{\tau_{i_1,j_1}\dots\tau_{i_k,j_k}}\\
x& \mapsto & x\tau_{i_1,j_1}\dots\tau_{i_k,j_k}
\end{array}
\end{displaymath}
is obviously a bijection, there exists $w'\in
\GH_{\tau_{1,2}\dots\tau_{2k-1,2k}}$ for which we have
$w=w'\tau_{i_1,j_1}\dots\tau_{i_k,j_k}$. From the definition of
$g$ we have   $g(n+2p)=g(n)+3p$ for all positive 
integers $n$ and $p$. Hence
\begin{displaymath}
\ls(w)\leq k+\ls(w')\leq k+(k-1)+g(n-2(k-1))=g(n)+2-k\leq g(n).
\end{displaymath}

{\bf Case~2:} $\lmod\{i_k,j_k\}\cap\{1,\dots,2k\}\rmod\leq 1$.
Without loss of generality we may also further assume  
$\{i_k,j_k\}\cap\{1,\dots,2k-2\}=\varnothing$.
Then, using Lemma~\ref{lm:aux}, we see that
there exists $\mu$ such that $\len(\mu)\leq 2(k-1)$ and
\begin{displaymath}
\tau_{1,2}\dots\tau_{2k-3,2k-2}\tau_{i_k,j_k}\mu
\GL\tau_{i_1,j_1}\dots\tau_{i_k,j_k}.
\end{displaymath}
This implies that the map
\begin{displaymath}
\begin{array}{rcl}
\GL_{\tau_{1,2}\dots\tau_{2k-3,2k-2}\tau_{i_k,j_k}} & \to & 
\GL_{\tau_{i_1,j_1}\dots\tau_{i_k,j_k}}\\
x& \mapsto & x\mu
\end{array}
\end{displaymath}
is a bijection. In particular, there exists
$v\in \GL_{\tau_{1,2}\dots\tau_{2k-3,2k-2}\tau_{i_k,j_k}}$
such that $v\mu=w$. Since $\tau_{1,2}\dots\tau_{2k-1,2k}w=w$,
it follows that $\tau_{1,2}\dots\tau_{2k-1,2k}v=v$ (because of
$\corank(\varphi(v))=\corank(\varphi(w))$). Hence
$v\in \GR_{\tau_{1,2}\dots\tau_{2k-1,2k}}$. 

From the above we derive $\tau_{2a-1,2a}v=v\tau_{2a-1,2a}=v$ for all $a=1,\dots,k-1$. Hence we can write $v=\tau_{1,2}\dots\tau_{2k-3,2k-2}v'$,
where $v'$ such that $\corank(\varphi(v'))=2$ and 
$\{m,m'\}\in \varphi(v')$ for all $m\leq 2k-2$. Using induction on
$n$ and the case of corank $2$ considered above, we obtain
$\ls(v')\leq g(n-2(k-1))$. Hence $\ls(v)\leq (k-1)+g(n-2(k-1))$
and from $w=v\mu$ we get 
\begin{displaymath}
\ls(w)\leq \ls(v)+\ls(\mu)\leq (k-1)+g(n-2(k-1))+2(k-1)=g(n).
\end{displaymath}
This completes the proof.
\end{proof}

%%%%%%%%%%%%%%%%%%%%%%%%%%%%%%%%%%%%%%%%%%%%%%%%%%%%%%%%%%%%%%%%%%%%%%%%%%%%%%%%%%%%%%%%%%%%%%%%%%%%%%%

\bigskip
\noindent
V.Mal.: Algebra, Department of Mechanics and Mathematics,
Kyiv Taras Shev\-chenko University,
Volodymyrska str., 64, UA-01033, Kyiv, UKRAINE,
e-mail: {\tt vmaltcev\symbol{64}univ.kiev.ua}

\bigskip
\noindent
V.Maz.: Department of Mathematics,
Uppsala University,
Box 480, SE-75106, Uppsala, SWEDEN,
e-mail: {\tt mazor\symbol{64}math.uu.se}


\begin{thebibliography}{9999}

\bibitem[Aiz]{Aiz} A. A\v\i zen\v stat, Defining Relations of Finite 
Symmetric Semigroups (in Russian), Mat. Sb. N. S., {\bf 45} (1958), 
261--280.

\bibitem[Br]{Br} R. Brauer, On algebras which are connected with the 
semisimple continuous groups, Ann. of Math. (2) {\bf 38} (1937), no. 4, 857--872.

\bibitem[CP]{CP} A.~Clifford, G.~Preston, The algebraic theory of 
semigroups. Vol. I. Mathematical Surveys, No. 7 American Mathematical 
Society, Providence, R.I. 1961

\bibitem[East]{East} J. East, A Presentation of the Singular part of the Symmetric Inverse Monoid, Preprint, 2004.

\bibitem[KM1]{KM1} G. Kudryavtseva, V. Mazorchuk, On conjugation in some transformation and Brauer-type semigroups, Preprint, Uppsala University, to appear in Publ. Math. Debrecen.

\bibitem[KM2]{KM2} G. Kudryavtseva, V. Mazorchuk, On presentation of Bra\-uer-\-type mo\-no\-ids, Preprint, Uppsala University, 2005, to appear
in Central Europ. J. Math.

\bibitem[KMM]{KMM} G. Kudryavtseva, V. Maltcev, V. Mazorchuk, $\GL$- and $\GR$-cross-sections in the Brauer semigroup, Preprint, Uppsala University, to appear in Semigroup Forum.

\bibitem[Mal]{Mal} V. Maltcev, Generating systems, ideals and the chief series of the Brauer semigroup, Proceedins of Kyiv University, Physical and Mathematical Sciences, 2004, no. 2, 59--65.

\bibitem[Maz1]{Maz1} V. Mazorchuk, On the structure of Brauer semigroup and its partial analogue, Problems in Algebra {\bf 13} (1998), 29--45.

\bibitem[Maz2]{Maz2} V. Mazorchuk, Endomorphisms of $\B_n$, $\PB_n$ and $\C_n$, Comm. Algebra {\bf 30} (2002), no. 7, 3489--3513.

\end{thebibliography}
\end{document}